\RequirePackage{fix-cm}
\documentclass[smallextended, envcountsect]{svjour3}       
\usepackage[utf8]{inputenc}
\usepackage[T1]{fontenc}
\smartqed  
\usepackage{graphicx,color,colortbl}
\usepackage{enumitem}
\usepackage{caption}
\usepackage{subcaption}
\usepackage{docmute} 
\renewcommand\thefootnote{*\arabic{footnote}}
\newcommand{\ovlU}{\overline{U}}
\newcommand{\tteal}{\textcolor{black}}
\newcommand{\tblue}{\textcolor{black}}

\newcommand{\bx}{\bar{x}}
\newcommand{\oLambda}{\overline{\Lambda}}
\newcommand{\blambda}{\bar{\lambda}}

\newcommand{\secondorder}{second-order }
\newcommand{\Lipschitz}{Lipschitz }

\newcommand{\Frobenius}{Frobenius }
\newcommand{\semidefinite }{semidefinite }

\newcommand{\worstcase}{worst-case }
\newcommand{\comple}{complementarity }
\newcommand{\statio}{stationarity }
\newcommand{\opti}{optimality }

\newcommand{\trust}{trust-region }
\newcommand{\primaldual}{primal-dual }
\newcommand{\STEPSIZE}{step size}
\definecolor{negative_curvature}{RGB}{187, 85, 102}
\usepackage{amsmath,amssymb,amsfonts,bm}
\usepackage{thmtools} 
\usepackage{siunitx, physics}
\usepackage{optidef}
\makeatletter
\let\c@proposition\c@theorem
\let\c@corollary\c@theorem
\let\c@lemma\c@theorem
\let\c@definition\c@theorem
\let\c@remark\c@theorem
\makeatother
\declaretheorem[
  name=Assumption,
  style=thmsty,
]{assumption}
\newlist{enuminasm}{enumerate}{1}
\setlist[enuminasm]{
  label=(\roman*),
  ref=\theassumption~(\roman*),
  leftmargin=0.08\textwidth,
  topsep=0.2\baselineskip,
  itemsep=0.2\baselineskip,
}
\newlist{enumlemma}{enumerate}{1}
\setlist[enumlemma]{
  label=(\arabic*),
  ref=\thelemma~(\arabic*),
  leftmargin=0.08\textwidth,
  topsep=0.2\baselineskip,
  itemsep=0.2\baselineskip,
  font=\normalfont
}
\newlist{enumprop}{enumerate}{1}
\setlist[enumprop]{
  label=(\arabic*),
  ref=\theproposition~(\arabic*),
  leftmargin=0.08\textwidth,
  topsep=0.2\baselineskip,
  itemsep=0.2\baselineskip,
  font=\normalfont
}

\DeclareMathOperator*{\argmin}{arg\,min}
\newcommand{\tangentspace}{T}
\newcommand{\heps}{\widehat{\varepsilon}}
\newcommand{\hnu}{\widehat{\nu}}
\newcommand{\ovzeta}{\overline{\zeta}}
\newcommand{\nnorm}[1]{\|#1 \|}
\newcommand{\fnorm}[1]{\norm{#1}_{\mathrm{F}}}
\newcommand{\paren}[1]{\qty(#1)}
\newcommand{\sbra}[1]{\qty{#1}}
\newcommand{\relmiddle}[1]{\mathrel{}\middle#1\mathrel{}}
\newcommand{\xk}{x^{k}}
\newcommand{\xell}{x^{\ell}}
\newcommand{\Xk}{X_{k}}

\newcommand{\Xkp}{{X_{k}^{P}}}
\newcommand{\invXkP}{\paren{X_{k}^{P}}^\dagger}
\newcommand{\Xkn}{{X_{k}^{N}}}
\newcommand{\invXkN}{\paren{X_{k}^{N}}^\dagger}

\newcommand{\Xx}{X}
\newcommand{\invXx}{X^{-1}}
\newcommand{\Zz}{Z}

\newcommand{\calNellx}{\mathcal{N}_1(\xell)}
\newcommand{\calNellZ}{\mathcal{N}_2(\Zell)}
\newcommand{\Lastk}{\mathcal{L}_k}
\newcommand{\calO}{\mathcal{O}}
\newcommand{\calD}{\mathcal{D}}
\newcommand{\invXk}{X^{-1}_{k}}

\newcommand{\invZk}{Z^{-1}_{k}}

\newcommand{\invXell}{X^{-1}_{\ell}}

\newcommand{\invXellp}{X^{-1}_{\ell+1}}

\newcommand{\rast}{r_{\ast}}
\newcommand{\xkp}{x^{k+1}}
\newcommand{\xellp}{x^{\ell+1}}
\newcommand{\Xell}{X_{\ell}}

\newcommand{\Zmu}{Z_{\mu}}

\newcommand{\Ker}{\mathop{\rm Ker}}
\newcommand{\Gast}{G_{\ast}}
\newcommand{\mineig}[1]{\lambda_{\min}\paren{#1}}
\newcommand{\maxeig}[1]{\lambda_{\max}\paren{#1}}
\newcommand{\eigi}[1]{\lambda_{i}\paren{#1}}
\newcommand{\xinit}{x^{1}}
\renewcommand{\order}[1]{{\rm O}\hspace{-0.2em}\left(#1\right)}

\newcommand{\innerstepsize}{\alpha_{\ell}}

\newcommand{\psimunu}[1]{\psi_{\mu,\nu} \paren{#1}}
\newcommand{\psimunutwo}[1]{\psi_{\mu,0} \paren{#1}}

\newcommand{\hesspsimu}[1]{\nabla^{2}_{xx} \psi_{\mu,\nu} \paren{#1}}
\newcommand{\hesspsimutwo}[1]{\nabla^{2}_{xx} \psi_{\mu,0} \paren{#1}}

\newcommand{\innerdirection}{d_{x^{\ell}}}

\newcommand{\innerZdirection}{d_{Z_{\ell}}}
\newcommand{\nablaij}[1]{\paren{\nabla^2_{xx} #1}_{ij}}
\newcommand{\Lagrange}[1]{L\paren{#1}}
\newcommand{\HessLagrange}[1]{\nabla^2_{xx} \Lagrange{#1}}
\renewcommand{\trace}[1]{\,\mathrm{tr}\paren{#1}}
\renewcommand{\det}[1]{\,\mathrm{det}\paren{#1}}

\newcommand{\invX}[1]{X(#1)^{-1}}
\newcommand{\uk}{u^k}
\newcommand{\ubar}{\bar{u}}

\newcommand{\quadf}[2]{{#2}^{\top} #1 #2}

\newcommand{\xl}{x^k}

\newcommand{\Dk}{{D}_{k}}
\newcommand{\Xast}{{X}_{\ast}}
\newcommand{\Xastdagger}{{X}_{\ast}^\dagger}
\newcommand{\Dast}{{D}_{\ast}}

\newcommand{\Vk}{{U}_{k}}
\newcommand{\Vast}{{U}_{\ast}}
\newcommand{\Dkp}{{D}_{k}^{P}}
\newcommand{\invDkp}{\paren{{D}_{k}^{P}}^{-1}}
\newcommand{\Dastp}{{D}_{\ast}^{P}}

\newcommand{\Vkp}{{U}_{k}^{P}}
\newcommand{\Vastp}{{U}_{\ast}^{P}}
\newcommand{\Dkz}{{D}_{k}^{N}}
\newcommand{\invDkz}{\paren{{D}_{k}^{N}}^{-1}}
\newcommand{\Dastz}{{D}_{\ast}^{N}}
\newcommand{\Vkz}{{U}_{k}^{N}}
\newcommand{\Vastz}{{U}_{\ast}^{N}}

\newcommand{\transpose}[1]{\paren{#1}^{\top}}

\newcommand{\muk}{\mu_k}
\newcommand{\mukp}{\mu_{k+1}}

\newcommand{\logdet}[1]{\mathrm{log}\,\mathrm{det}#1}

\newcommand{\Zk}{Z_k}
\newcommand{\Zell}{Z_\ell}

\newcommand{\Zellp}{Z_{\ell+1}}

\newcommand{\q}[1]{\pdv{#1}{x_q}}
\newcommand{\euclid}{{\mathbb{R}}}
\newcommand{\nonnegative}{\mathbb{R}_{+}}

\newcommand{\symmetric}{\mathbb{S}}
\newcommand{\xast}{x^\ast}
\newcommand{\Acali}{\mathcal{A}_i}
\newcommand{\Acalj}{\mathcal{A}_j}
\newcommand{\Aast}{\mathcal{A}^\ast}
\newcommand{\Dx}[2]{\Delta X(#1;#2)}
\newcommand{\mupn}{\paren{1+\nu}\mu}
\newcommand{\mukpn}{\paren{1+\hnu(\muk)} \mu_k}
\newcommand{\feasi}{\mathcal{X}_{++}}
\newcommand{\convfeasi}{\mathrm{conv}\paren{\feasi}}
\newcommand{\calK}{\mathcal{K}}

\newcommand{\calKo}{\mathcal{K}}

\newcommand{\calL}{\mathcal{L}}

\newcommand{\epsmu}{\varepsilon_\mu}

\newcommand{\Proj}{\mathrm{Proj}}
\newcommand{\diag}{\mathrm{diag}}
\newcommand{\FPB}{F_{\mathrm{PB}}^\mu}
\newcommand{\FBC}{F_{\mathrm{BC}}^\mu}

\newcommand{\Hell}{H_\ell}

\newcommand{\hmin}{h_{\min}}
\newcommand{\hmax}{h_{\max}}
\newcommand{\calHell}[1]{\mathcal{H}_\ell\paren{#1}}
\newcommand{\iotamin}{\kappa_{\min}}
\newcommand{\iotamax}{\kappa_{\max}}

\newcommand{\pmunu}{\psi_{\mu, \nu}}
\newcommand{\xbar}{\bar{x}}
\newcommand{\Zbar}{\bar{Z}}

\usepackage[ruled, noline, linesnumbered,inoutnumbered, commentsnumbered, longend,procnumbered,norelsize]{algorithm2e}
\makeatletter
\let\original@algocf@latexcaption\algocf@latexcaption
\long\def\algocf@latexcaption#1[#2]{%
  \@ifundefined{NR@gettitle}{%
    \def\@currentlabelname{#2}%
  }{%
    \NR@gettitle{#2}%
  }%
  \original@algocf@latexcaption{#1}[{#2}]%
}
\makeatother
\SetAlgoProcName{Update\hspace*{-0.05cm}}{Update}
\SetAlgorithmName{Algorithm\hspace*{-0.05cm}}{Algorithm}{Algorithms}
\usepackage{cite}
\PassOptionsToPackage{hyphens}{url}
\usepackage{hyperref}
\usepackage{xcolor}
\definecolor{Gray}{gray}{0.9}
\hypersetup{%
  breaklinks = true,
  bookmarksnumbered = true,%
  hidelinks,%
  colorlinks = true,
  filecolor = magenta,%
  setpagesize = false,%
  colorlinks,
  linkcolor={black},
  citecolor={green!35!black},
  urlcolor={blue!50!black}
}

\makeatletter
\let \c@definition\c@proposition
\let \c@proposition\c@definition
\let \c@theorem\c@remark
\let \c@theorem\c@proposition
\let \c@proposition\c@lemma
\let \c@lemma\c@theorem
\makeatother
\usepackage{cleveref}
\crefname{assumption}{Assumption}{Assumptions}
\crefname{lemma}{Lemma}{Lemmas}
\Crefname{lemma}{Lemma}{Lemmas}
\crefalias{enuminasmi}{assumption}
\crefalias{enumlemmai}{lemma}
\crefname{proposition}{Proposition}{Propositions}
\Crefname{proposition}{Proposition}{Propositions}
\crefalias{enumpropi}{proposition}

\crefname{definition}{Definition}{Definitions}
\crefname{lemmabreak}{Lemma}{Lemmas}
\crefname{theorem}{Theorem}{Theorems}
\Crefname{theorem}{Theorem}{Theorems}
\crefname{remark}{Remark}{Remarks}
\Crefname{remark}{Remark}{Remarks}
\crefname{equation}{}{}
\Crefname{equation}{}{}
\crefname{appendix}{Appendix}{Appendices}
\Crefname{appendix}{Appendix}{Appendices}
\crefname{procedure}{Update}{Updates}
\Crefname{procedure}{Update}{Updates}
\crefname{algorithm}{Algorithm}{Algorithms}
\Crefname{algorithm}{Algorithm}{Algorithms}
\crefname{section}{Section}{Sections}
\Crefname{section}{Section}{Sections}
\newcommand{\RN}[1]{%
  \textup{\uppercase\expandafter{\romannumeral#1}}%
}
\spnewtheorem*{xproof}{}{\itshape}{\rmfamily}

\renewenvironment{proof}[1][\proofname]
{\renewcommand\xproofname{#1}\xproof}
{\endxproof}

\usepackage{color}
\def\ASOSP{approx.~SOSP}
\def\wSOSP{\rm{w}-SOSP}
\definecolor{darkgreen}{rgb}{0,0.5,0}
\definecolor{purple}{rgb}{1,0,1}
\usepackage{siunitx}

\begin{document}
\title{Complexity analysis of interior-point methods for second-order stationary points of nonlinear semidefinite optimization problems
  \footnote{This work was supported by the Japan Society for
  the Promotion of Science KAKENHI Grant Number 19H04069, 20K19748,
  20H04145, and 23H03351.  It was conducted when the first author was
  a student of the University of Tokyo and is irrelevant to his
  present affiliation.  } }
\author{Shun Arahata \and Takayuki Okuno \and
  Akiko Takeda}
\titlerunning{Interior-point methods for SOSPs of NSDPs}
\authorrunning{S. Arahata, T. Okuno, and A. Takeda}

\institute{ S. Arahata \at Graduate School of Information Science and Technology, The University of Tokyo, Tokyo, Japan\\
\email{sarahata@alumni.u-tokyo.ac.jp}
    \and
    T. Okuno
     \at Faculty of Science and Technology Department of Science and Technology, Seikei University, Tokyo, Japan
    \at Center for Advanced Intelligence Project, RIKEN, Tokyo, Japan\\
    \email{takayuki-okuno@st.seikei.ac.jp}
    \and
    A. Takeda \at Graduate School of Information Science and Technology, The University of Tokyo, Tokyo, Japan \at Center for Advanced Intelligence Project, RIKEN, Tokyo, Japan\\
    \email{takeda@mist.i.u-tokyo.ac.jp, akiko.takeda@riken.jp}
}

\date{Received: date / Accepted: date}

\maketitle
\begin{abstract}
    We propose a primal-dual interior-point method \allowbreak (IPM)  with convergence to second-order stationary points (SOSPs) of nonlinear semidefinite optimization problems, abbreviated as NSDPs.
    As far as we know, the current algorithms for NSDPs only ensure convergence to first-order stationary points such as Karush-Kuhn-Tucker points, but without a worst-case iteration complexity.
    The proposed method generates a sequence approximating SOSPs while minimizing a primal-dual merit function for NSDPs by using
    scaled gradient directions and directions of negative curvature.
    Under some assumptions, the generated sequence accumulates at an SOSP with a worst-case iteration complexity.
    This result is also obtained for a primal IPM with a slight modification.
    Finally, our numerical experiments show the benefits of using directions of negative curvature in the proposed method.
\end{abstract}

\section{Introduction}
\label{intro}

We consider the following nonlinear \semidefinite optimization problems (NSDPs), which are possibly nonconvex:
\begin{mini}
  {x\in \mathbb{R}^n}{f(x)}
  {\label{eq:main}}{}
  \addConstraint{X(x)}{\in \symmetric^m_{+},}
\end{mini}
where $f:\euclid^n \to \euclid$ and $X:\euclid^n \to \symmetric^m$ are twice continuously differentiable functions and
$\euclid^n$ and $\symmetric^m$ denote
the spaces of $n$-dimensional real vectors and $m\times m$ real symmetric matrices, respectively. Moreover,
$\symmetric^m_{++(+)}$ denotes the set of
positive (semi)definite matrices in $\symmetric^m$.
By restricting the range of $X$ onto the space of diagonal matrices,
NSDP\,\eqref{eq:main} reduces to the standard inequality-constrained nonlinear optimization problem (NLP).
When all the functions of a semidefinite optimization problem are affine, the problem is a linear semidefinite optimization problem (LSDP) and has been studied extensively~\cite{wolkowicz2012handbook}.
\par
Though studies of NSDPs are still fewer than those of LSDPs and NLPs, they are of great
  importance from a practical point of view.
Indeed, NSDPs arise from various application fields including control~\cite{Bassem2001,Hoi2003,Kocvara2005}, statistics~\cite{Houduo2006}, finance~\cite{Konno2003,KonnoKawadai2003,Leibfritz2008}, and structural optimization~\cite{Bendsoe1994,Michal2004,Kanno2006}.
Positive semidefinite matrix factorization problems~\cite{Vandaele2018,Lahat2020} and rank minimization problems~\cite{Fazel2003} are also important applications of NSDPs.
Moreover, NSDPs have been studied in terms of optimality conditions: For example, the Karush-Kuhn-Tucker (KKT) conditions and the second-order conditions for NSDPs were studied in detail by Shapiro~\cite{shapiro1997first} and Forsgren~\cite{forsgren2000optimality}.
Further examples are: the strong second-order conditions by Sun~\cite{sun2006strong}, the local duality by Qi~\cite{qi2009local}, sequential optimality conditions by Andreani et al.~\cite{andreani2020optimality}, and the reformulation of optimality conditions via slack variables by Louren{\c{c}}o et al.~\cite{Lourenco2018}. We also refer to the book by Bonnans and Shapiro~\cite{BonnansShapiro2000}.
Supported by those theories, various algorithms have been developed for NSDPs, including
primal interior-point methods (IPMs)~\cite{Jarre2000,Leibfritz2002}, primal-dual IPMs~\cite{Yamashita2012,Yamakawa2014,Yamashita2020,Okuno2020,Okuno2020112784,okuno2018primal}, augmented Lagrangian methods~\cite{Fukuda2018,andreani2020optimality}, and sequential quadratic semidefinite programming methods~\cite{zhao2016superlinear,freund2007nonlinear,correa2004global,Yamakawa2020}.
\par
Before the present study, the existing algorithms for NSDPs only ensure global convergence to first-order stationary points\footnote{
{After the present study was submitted to the arXiv as a preprint\,\cite{arahata2021interior}, Yamashita\,\cite{yamashita2022convergence} presented another interior-point method designed in the framework of the trust region method with convergence to an SOSP but without complexity analysis. }}
such as KKT points, AKKT points, and TAKKT points. (The AKKT and TAKKT points are new optimality concepts presented most recently in \cite{andreani2020optimality}.)
In contrast, many algorithms for computing
second-order stationary points, SOSPs in short, of NLPs have been proposed for several decades, e.g., for unconstrained problems, the negative curvature method by
McCormick~\cite{McCormick1977}, the trust-region method by Sorensen~\cite{Sorensen1982}, and the cubic regularized Newton method by Nesterov and Polyak~\cite{Nesterov2006}. Even for constrained problems, we can find penalty methods by Auslender~\cite{Auslender1979} and by Facchinei and Lucidi~\cite{FacchineiLucidi1998}, the squared slack variables technique by Mukai and Polak~\cite{Mukai1978}, and the negative curvature method with projection by Goldfarb et al.~\cite{Goldfarb2017}.
\par
Recently, motivated by applications in machine learning, quite a few algorithms for NLPs with worst-case iteration complexities for an approximate SOSP are proposed.
For example, with respect to unconstrained optimization,  we find \cite{Curtis2017,CurtLubbRobi18,curtis2020trustregion} using the \trust method and \cite{Curtis2019} using the negative curvature method.
For constrained optimization, there also exist many algorithms depending on constraints. \Cref{table:constrained} summarizes the algorithms for constrained NLPs with convergence to SOSPs.
\par
One may think of transforming NSDPs to NLPs  by means of squared slack variables.
Indeed, the semidefinite constraint $X(x)\in \symmetric^m_{+}$ can be equivalently transformed
into the equality constraint $X(x)=Y^2$ with a slack variable $Y\in \symmetric^m$.
Hence, NSDP~\eqref{eq:main} reduces to the equality-constrained NLP possessing $(x,Y)$ as variables,
and thus existing NLP algorithms with convergence to SOSPs are applicable to the transformed NLP.
However, such an approach has drawbacks, e.g., the number of variables increases and, as indicated in \cite{Lourenco2018},
there may exist a discrepancy between the set of SOSPs of the NSDP and that of SOSPs of the transformed NLP.
\begin{table}[htbp]
  \centering
  \caption{Constrained NLP algorithms with convergence to SOSPs
({\small The terms ``eq'' and ``ineq'' in the Constraint-column stand for equality and inequality,~respectively.})}
  \begin{tabular}{lll}
    Constraint          & Algorithm                                                                        & Complexity \\ \hline
    Nonlinear eq        & Penalty~\cite{El-Alem1996} and projection~\cite{Goldfarb2017}                    & No         \\
    Linear ineq         & Active-set~\cite{Foresgren1997}                                                  & No         \\
    Nonlinear ineq      & Squared slack~\cite{Mukai1978}, penalty~\cite{Auslender1979,FacchineiLucidi1998} & No         \\
                        & augmented Lagrangian~\cite{DiPillo2005, Andreani2010}                            &            \\
                        & {\bf primal-dual IPM\cite{Conn2000,Moguerza2003}}                                &            \\
    Nonlinear eq        & Projection~\cite{SunFazel2018}                                                   & Yes        \\
                        & and proximal augmented Lagrangian~\cite{Xie2020}                                 &            \\
    Linear ineq         & Active-set~\cite{Lu2019} and trust-region~\cite{NouiehedRazaviyayn2020}          & Yes        \\
    Closed convex set   & Second-order Frank-Wolfe~\cite{Mokhtari2018,Nouiehed2020,Hallak2020}                        & Yes        \\
    Nonnegative orthant & {\bf Primal IPM~\cite{ONeilWright2020}}                                          & Yes        \\
    Nonlinear ineq      & {\bf Primal IPM~\cite{Hinder2020}}                                               & Yes
  \end{tabular}
  \label{table:constrained}
\end{table}
\subsection{Contributions}
The main contributions in this paper are summarized in the following two theoretical points.
\begin{enumerate}[label = \arabic*)]
  \item {\it SOSPs of NSDPs}\hspace{7pt}
        We present the first primal-dual IPM with convergence to SOSPs of NSDP~\eqref{eq:main}.
        In the primal-dual strictly feasible region $\{(x,Y)\mid X(x)\in \symmetric^m_{++}, Y\in \symmetric^m_{++}\}$, the primal-dual IPM generates a sequence of approximate SOSPs, abbreviated as {\ASOSP}s, which accumulate at an SOSP of NSDPs under some assumptions. Each {\ASOSP} is defined using the primal-dual merit function given in~\cite{Yamashita2012}. In order to compute such {\ASOSP}s,
        we utilize scaled gradient directions\footnote{A scaled gradient direction is intended to be the steepest-descent one premultiplied with a positive definite symmetric matrix.} and directions of negative curvature.
        The difficulty of theoretical analysis of NSDPs often stems from the so-called sigma term (see Definition\,\ref{def:SONC}), which reflects the curvature of the positive semidefinite cone. We need to handle it carefully in order to show the convergence to SOSPs of NSDPs.
  \item {\it a worst-case iteration complexity for {\ASOSP}s}\hspace{7pt}
        We give a \worstcase iteration complexity of the proposed \primaldual IPM for
        computing an {\ASOSP}.
        Our {\STEPSIZE} rule makes it possible to keep the generated iterates in the strictly feasible region and obtain a \worstcase iteration complexity.
        Another point to note is that this result is established in the primal-dual framework, whereas previous studies~\cite{ONeilWright2020,Hinder2020} regarding NLPs gave \worstcase iteration complexities for primal IPMs.
\end{enumerate}
We stress that the proposed primal-dual IPM is the first NSDP algorithm for SOSP equipped with a worst-case iteration complexity.
We can obtain a primal IPM with convergence to SOSPs and a worst-case iteration complexity with a slight modification.
\subsection{Related works}
\paragraph{Primal-dual IPM for NSDPs}
Let us review some existing studies on primal-dual IPMs (PDIPMs in short) for NSDPs as the most relevant works.
Similar to the (path-following) PDIPM for LSDPs (see e.g.,\cite[Chapter~10]{wolkowicz2012handbook}), the fundamental framework of the existing PDIPMs for NSDPs is to approach a KKT triplet of \tteal{an NSDP by computing} perturbed KKT triplets and driving a perturbation parameter to zero.
We believe that Yamashita, Yabe, and Harada~\cite{Yamashita2012} presented the first PDIPM for NSDPs and showed its global convergence to a KKT triplet of an NSDP using the family of Monteiro-Zhang directions. Its local convergence property was analyzed by Yamashita and Yabe in \cite{yamashita2012local}.
Afterward, Kato et al.~\cite{kato2015interior} studied the global convergence of PDIPMs using a different penalty function from \cite{Yamashita2012}.
Yamakawa and Yamashita~\cite{Yamakawa2014} developed a different PDIPM based on the shifted barrier KKT conditions for NSDPs.
Okuno~\cite{Okuno2020} analyzed \tteal{the} local convergence of a PDIPM using the family of Monteiro-Tsuchiya directions.
{ The trust-region-based PDIPM was presented by Yamashita
  et al.~\cite{Yamashita2020}. It is worth mentioning that no NSDP
  algorithms including IPMs, except for linear SDPs, \tteal{are} equipped with
  iteration-complexity results.  }
\paragraph{IPM with SOSPs and iteration complexities for NLPs}
We review some existing IPMs with convergence to an SOSP of constrained NLPs.
The fundamental idea of such IPMs is to approach the set of SOSPs by tracking a path-like set formed by SOSPs of the reformulated problems obtained by means of interior penalty functions such as a log-barrier function.
To compute SOSPs of such penalized NLPs, Conn et al.~\cite{Conn2000} used the trust-region method in the framework of PDIPM and 
Moguerza and Prieto~\cite{Moguerza2003} employed the negative curvature method coupled with the modified Newton method.
Recently, worst-case iteration complexities of primal IPMs to SOSPs were analyzed by O'Neill and Wright~\cite{ONeilWright2020} for NLPs with the nonnegative orthant constraint $x\ge 0$ using the Newton-CG method~\cite{Royer2020}, and also by Hinder and Ye~\cite{Hinder2020} for NLPs with general nonlinear inequalities using the trust-region method.
We note that there are still no works handling worst-case iteration complexities for PDIPMs.
The above IPM papers are shown in bold in the aforementioned~\Cref{table:constrained}.
\subsection{Notation and outline of the paper}
\label{subsection:notation}
Throughout this paper, the following symbols are often used.
For $A\in \symmetric^m$, $\lambda_i(A)$ denotes its $i$-th largest eigenvalue. In particular, $\mineig{A}$ and $\maxeig{A}$ denote the minimum and maximum eigenvalues of $A$, respectively. For $B \in \symmetric^m_{+}$, $B^{\frac{1}{2}} $ denotes the positive \semidefinite root of $B$, that is, $B^{\frac{1}{2}} \in \symmetric^m_{+}$ and $\paren{B^{\frac{1}{2}}}^2 = B$.
For a matrix $C\in \euclid^{m\times n}$, the transpose, the Moore-Penrose pseudo-inverse, the rank, and the null or kernel space of $C$ are denoted by $C^\top, C^\dagger, \rank (C),$ and ${\rm Ker}(C)$, respectively.
Also, $\fnorm{C}$ denotes the \Frobenius norm of $C$, that is, $\fnorm{C} = \sqrt{\sum_{i=1}^m \sum_{j=1}^n C_{ij}^2}$, where $C_{ij}$ denotes the $(i,j)$-th entry of $C$. When $C$ is a $m \times m$ real symmetric matrix, $\fnorm{C}$ equals to $\sqrt{\sum_{i=1}^m \lambda_i(C)^2}$.
In addition, $\norm{C}$ denotes the spectral norm, that is, $\norm{C} = \sqrt{\lambda_{\max}(C^\top C)}$. For a vector $c \in \euclid^n$, $\norm{c}$ denotes the Euclidean norm of $c$, and $c_i$ denotes the $i$-th element of the vector. For matrices $Y\in \euclid^{ m \times n }$ and $Z\in \euclid^{m \times n}$, their inner product is denoted as $\langle Y, Z \rangle \coloneqq \trace{Y^\top Z}$.
\par
We also define $\mathcal{A}_i (x) \coloneqq  \pdv{X}{x_i}\paren{x}$ for each $i \in \sbra{1,\ldots, n}$. In addition, we define $\Aast (x)$ for $x \in \euclid^n$ as the following operator from $\symmetric^m$ to $\euclid^n$:
$$
\Aast(x) Z \coloneqq	
  \paren{
    \langle \mathcal{A}_i(x), Z\rangle, \ldots,
    \langle \mathcal{A}_n(x), Z \rangle
  }^\top,
$$
where $Z\in \symmetric^m$. Let $d\in \mathbb{R}^n$ and define
$
  \Delta X(x;d) \coloneqq \sum_{i=1}^n \mathcal{A}_i \paren{x} d_i
  .$
\par
\textbf{Outline.}
  The paper is organized as follows. The rest of this section introduces assumptions of  Lipschitz continuities. \Cref{section:opticond} reviews optimality conditions and constraint qualifications of NSDPs. Then, \Cref{chapter:algo} describes the proposed IPM and \Cref{chapter:convergence_outer,chapter:inner} analyze its convergence, while \Cref{sec:primal} states that, with a slight modification, it is possible to make a primal IPM with almost the same convergence properties. \Cref{section:numexp} presents the result of numerical experiments. \Cref{section:conclu} concludes the paper.
\subsection{Assumptions on the objective and constraint}
In order to establish a \worstcase iteration complexity for the proposed IPM, we make assumptions on \Lipschitz continuities and boundedness of derivatives. Let $\feasi$ denote the strictly feasible region, that is,
\begin{equation}
  \feasi \coloneqq \bigl\{x\in \euclid^n\;: \allowbreak X(x) \in \symmetric^m_{++}\bigr\}\label{eq:feasdef}
\end{equation}
and $\convfeasi$ denote the convex hull of $\feasi$. The following~\Cref{assumption:Lipschtiz} for the objective and constraint on $\convfeasi$ are implicitly supposed throughout the paper.
Recall that twice \tteal{continuous} differentiability of the objective and that of the constraint were assumed at the beginning of the paper.
\begin{assumption}\label{assumption:Lipschtiz}
  The gradient of $f$ is $\widetilde{L}_1$-\Lipschitz continuous on $\convfeasi$ and the Hessian is $\widetilde{L}_2$-\Lipschitz continuous on $\convfeasi$, that is, we have
  \begin{align}
    \norm{\nabla f(x) - \nabla f(z)}     & \leq\widetilde{L}_1 \norm{x-z} \label{eq:gradf},  \\
    \norm{\nabla^2 f(x) - \nabla^2 f(z)} & \leq \widetilde{L}_2 \norm{x-z} \label{eq:hessf},
  \end{align}
  for any $x,z \in  \convfeasi$.
  Moreover, we assume that, for any $x \in \convfeasi$,
  \begin{align}
    \sum_{i=1}^n \fnorm{\Acali \paren{x}}                        & \leq L_0 \label{eq:pdvXxiLips},            \\
    \sum_{i=1}^n \sum_{j=1}^n \fnorm{\pdv{X}{x_i}{x_j}\paren{x}} & \leq \widehat{L}_1 \label{eq:hessbounded},
  \end{align}
  where $\Acali(\cdot)$ is defined in \Cref{subsection:notation}.
  For the second-order derivative of $X$, we also assume that
  \begin{equation}
    \sum_{i=1}^n \sum_{j=1}^n \fnorm{\pdv{X}{x_i}{x_j}\paren{x}-\pdv{X}{x_i}{x_j}\paren{z}}
    \leq \widehat{L}_2 \norm{x-z}
    \label{as:hessXlips}
  \end{equation}
  for any $x,z \in \convfeasi$.
  For simplicity, we define
  $$L_1 \coloneqq \max \sbra{\widetilde{L}_1, \widehat{L}_1},\ L_2 \coloneqq \max \sbra{\widetilde{L}_2, \widehat{L}_2}.$$
\end{assumption}
The assumptions\,\eqref{eq:gradf} and \eqref{eq:hessf} on the objective function $f$ are often assumed in deriving iteration complexities~\cite{Goldfarb2017,Hinder2020}.
The constants in the remaining assumptions on the constraint function $X$ are ensured to exist if $X$ is affine or $\feasi$ is bounded.
\par
The boundedness of the derivatives implies the \Lipschitz continuities.
Indeed, under~\cref{assumption:Lipschtiz}, we can show
\begin{align}
  \fnorm{X(x) - X(z)}                                    & \leq L_0\norm{x - z} \label{eq:Xlip},      \\
  \sum_{i=1}^n \fnorm{\Acali \paren{x}-\Acali \paren{z}} & \leq L_1 \norm{x -z}, \label{eq:gradxlips}
\end{align}
for any $x,z \in \convfeasi$.

\section{Optimality conditions for NSDPs}
\label{section:opticond}
This section aims to introduce \opti conditions together with constraint qualifications for NSDP\,\eqref{eq:main}.
First, we define the Fritz-John conditions.
\begin{definition}[Fritz-John (FJ) conditions]
    Let
    $L^g(x, \lambda, \Omega) \coloneqq \lambda f(x) - \langle X(x) , \Omega \rangle$
    be the generalized Lagrangian of NSDP~\eqref{eq:main}.
    We say that the FJ conditions for NSDP\,\eqref{eq:main} hold
    at $\bar{x}\in \euclid^n$ if there exists some $\paren{\bar{\lambda}, \overline{\Omega}} \in \euclid \times \symmetric^m$ such that
    \begin{subequations}
        \label{def:FJ}
        \begin{align}
            \nabla_{x} L^g(\bar{x}, \bar{\lambda}, \overline{\Omega}) & = 0, \label{eq:FJstationary}                                 \\
            X(\bar{x})\overline{\Omega}                               & = O, \label{eq:FJcomp}                                       \\
            X(\bar{x})                                                & \in \symmetric^m_{+}, \label{eq:FJfeasi}                     \\
            (\bar{\lambda}, \overline{\Omega})                        & \in \nonnegative \times \symmetric^m_{+}, \label{eq:FJdualf} \\
            \bar{\lambda} + \fnorm{\overline{\Omega}}                 & \neq 0,\label{eq:FJnonzero}
        \end{align}
    \end{subequations}
    where $\nabla_{x} L^g(\bar{x}, \bar{\lambda}, \overline{\Omega}) =
    \overline{\lambda} \nabla f(\bx) - \Aast\paren{\bx} \overline{\Omega}.$ Here $\Aast (\cdot)$ is defined in \Cref{subsection:notation}.
We call $\bar{x}$ an FJ point for NSDP\,\eqref{eq:main}.
\end{definition}
The FJ conditions are necessary conditions for the local optimality regardless of any constraint qualifications~\cite[Proposition 5.87]{BonnansShapiro2000}.
Next, we introduce the Mangasarian–Fromovitz constraint qualification (MFCQ)~\cite[§ 5.3.4]{BonnansShapiro2000} to describe the KKT conditions for NSDP\,\eqref{eq:main}.
\begin{definition}[Mangasarian–Fromovitz constraint qualification]
    Let $x \in \euclid^n$ be a feasible point for NSDP\,\eqref{eq:main}.
    If there exists a vector $h\in \mathbb{R}^n$ such that $X(x) + \Delta X(x; h) \in \symmetric^m_{++},$ where $\Delta X(x;h)$ is defined in \cref{subsection:notation}, then we say that the MFCQ holds at $x$.
\end{definition}
The KKT conditions for NSDP\,\eqref{eq:main} can be formally defined as follows.
\begin{definition}[Karush–Kuhn–Tucker (KKT) conditions]
    Define the Lagrangian of NSDP\,\eqref{eq:main} by $L(x, \Lambda) \coloneqq f(x) - \langle X(x), \Lambda \rangle$ for $(x,\Lambda)\in \euclid^n\times \symmetric^m$.
    We say that the KKT conditions for NSDP\,\eqref{eq:main} hold at $\bx$ if there exists some $\oLambda \in \symmetric^m$ satisfying
    \begin{subequations}\label{eq:KKT}
        \begin{align}
            \nabla_x L(\bx, \oLambda)  = 0 & \;\text{\normalfont(stationarity of the Lagrangian)},\label{eq:stationary} \\
            X(\bx) \in \symmetric^m_{+}    & \; \text{\normalfont(primal feasibility)},\label{eq:prfeasible}            \\
            \oLambda \in \symmetric^m_{+}  & \; \text{\normalfont(dual feasibility)},\label{eq:dufeasible}              \\
            X(\bx) \oLambda= O             & \; \text{\normalfont(complementarity)},\label{eq:complementarity}
        \end{align}
    \end{subequations}
    where $\nabla_x L(\bx, \oLambda) = \nabla f(\bx)
        -\Aast\paren{\bx}\oLambda.$ We call $\bx$ and $\oLambda$ a
    KKT point and a Lagrange multiplier matrix, respectively.
\end{definition}
Hereafter, for a KKT point $\bar{x}$, we denote by $\Lambda(\bar{x})$
the set of Lagrange multiplier matrices satisfying the KKT conditions.
It is well known that, under constraint qualifications such as the MFCQ, $\Lambda(\bar{x})\neq \emptyset$ holds at
a local optimum $\bar{x}$ of NSDP\,\eqref{eq:main}, 
namely, $\bar{x}$ is a KKT point.
See \cite[Theorem 5.84]{BonnansShapiro2000} for the proof.
It is worth mentioning that the MFCQ implies that
$\Lambda(\bar{x})$ is nonempty and bounded.
\par
Next, we state the second-order necessary conditions and the ``weak'' second-order necessary conditions.
\begin{definition}[second-order necessary conditions]\label{def:SONC}
  Let $\bx\in \euclid^n$ be a KKT point for NSDP\,\cref{eq:main} { and $r \coloneqq \rank (X(\bx))$. Take a matrix
  $\ovlU\in \euclid^{m\times (m-r)}$ such that the columns form an orthonormal basis of $\Ker(X(\bx))$.}
  We say that $\bx$ satisfies the second-order necessary conditions if
    \begin{equation}
        \sup_{\Lambda \in \Lambda (\bx)} d^\top \paren{\nabla^2_{xx} L(\bx,\Lambda) + H(\bx,\Lambda)}d \geq 0,\ \forall d \in C(\bx)
        \label{def:SOSPwithoutSC}
    \end{equation}
    holds and then call such $\bx$ an {SOSP} of NSDP\,\eqref{eq:main},
    where 
$H(\bx,\Lambda) \in \symmetric^m $ is defined as
    \begin{equation}
        \paren{H(\bx,\Lambda)}_{ij} \coloneqq  2\trace{\Acali(\bx) X(\bx)^\dagger \Acalj(\bx) \Lambda}
        \label{eq:sigma-term}
    \end{equation}
    for each $i,j$ and $X(\bx)^\dagger$ means the Moore-Penrose pseudo-inverse of $X(\bx)$. The expression $d^\top H(\bx,\Lambda) d$ is referred to as a {\it sigma term}. 
Moreover, $C(\bx)$ denotes the critical cone at $\bx$, which is represented as 
\begin{align}
        C(\bx)                                         & = \sbra{d \in \euclid^n \relmiddle| \Dx{\bx}{d} \in T_{\symmetric^m_{+}}(X(\bx)), \nabla f(\bx)^\top d = 0},\notag 
        \end{align}
where $T_{\symmetric^m_{+}}(X(\bx))$ denotes the tangent cone of 
        $\symmetric^m_{+}$ at $X(\bx)$ and is represented as\footnote{
When $\Ker(X(\bx))=\{0\}$, thus in the case where 
$X(\bx)\in \symmetric^m_{++}$, 
we write $\tangentspace_{\symmetric^m_{+}}\paren{X(\bx)}=
\symmetric^m$.} 
        \begin{align}
        \tangentspace_{\symmetric^m_{+}}\paren{X(\bx)} & = \sbra{E \in \symmetric^m \relmiddle| \ovlU^\top E \ovlU \in \symmetric^{m-r}_{+}}.\notag
    \end{align}
\end{definition}
{We refer readers to \cite[§5.3.1]{BonnansShapiro2000} for more details of the sigma term, critical cone, and tangent cone for the NSDP.}
Let $\bar{x}$ be a local optimum of NSDP\,\eqref{eq:main}. If the MFCQ holds at $\bar{x}$, then $\bar{x}$ is an SOSP~\cite[Theorem 5.88]{BonnansShapiro2000}.
\par
Next, we introduce a weaker concept than SOSP.
\begin{definition}[weak second-order necessary conditions]
    Let $\bx\in \euclid^n$ be a KKT point for NSDP\,\cref{eq:main}.
    We say \allowbreak that the weak second-order necessary conditions hold at $x$ if inequality~\eqref{def:SOSPwithoutSC}
    holds for all $d \in \mathcal{L}(\bx)$, where $\mathcal{L}(\bx)$ is the linear subspace defined by
\begin{equation}
        \mathcal{L}(\bx) \coloneqq
        \{ d \in \mathbb{R}^n \mid  \ovlU^{\top} (\Dx{\bx}{d}) \ovlU = O
        \}\label{eq:linearspaceL}
    \end{equation}
    if $\Ker(X(\bx))\neq \{0\}$, where $\ovlU$ is the same matrix as in the definition of  $\tangentspace_{\symmetric^m_{+}}\paren{X(\bx)}$. In particular, we define $\mathcal{L}(\bx)\coloneqq \euclid^n$ 
if $\Ker(X(\bx))=\{0\}$. 
We call $\bx$ a {\wSOSP} of NSDP\,\eqref{eq:main}.
\end{definition}
Note that 
when $X(\bx)$ has at least one zero-eigenvalue, that is,
$\bx$ is situated on the topological boundary of $\feasi$, the condition
$\Ker(X(\bx))\neq \{0\}$ holds.
\par
An SOSP is actually a {\wSOSP}. Indeed, let $\bx$ be an SOSP.
Choose $\Lambda\in \Lambda(\bx)$ arbitrarily, and let $\ovlU=[\ovlU_1,\ovlU_2]$ be a partition of $\ovlU$ 
such that $\Lambda$ is represented as $\Lambda = \ovlU_1\Theta\ovlU_1^{\top}$ with a positive definite matrix $\Theta$.\footnote{
If $\Lambda\neq O$, such $\ovlU_1$ exists.}  According to \cite[equation\,(5.219)]{BonnansShapiro2000}, it follows that 
\begin{equation}
\begin{split}
C(\bx)=\{
d\in \euclid^n\mid 
\ovlU_1^{\top}\Dx{\bx}{d}\ovlU_1=O,\ \ovlU_1^{\top}&\Dx{\bx}{d}\ovlU_2=O,\\
&\ovlU_2^{\top}\Dx{\bx}{d}\ovlU_2\in \symmetric_{+}^{\rho}
\},
\end{split}
\label{eq:cxu1u2}
\end{equation}
where $\rho$ is the number of columns of $\ovlU_2$.
Notice that 
$\mathcal{L}(\bx)$ is obtained by replacing 
$\ovlU_2^{\top}\Dx{\bx}{d}\ovlU_2\in \symmetric_{+}^{\rho}$ with $\ovlU_2^{\top}\Dx{\bx}{d}\ovlU_2=O$ in \eqref{eq:cxu1u2}, and hence from $\{O\}\subseteq \symmetric_{+}^{\rho}$ we see 
$\mathcal{L}(\bx)\subseteq C(\bx)$. Thus, any SOSP is a {\wSOSP}. On the other hand, the converse is not necessarily true.
\par
This gap, however, is completely filled 
under the strict complementarity condition at $\bx$, namely, there exists some $\oLambda\in \Lambda(\bar{x})$ such that $X(\bar{x})+\oLambda\in \symmetric^m_{++}$. 
This is because 
by setting $\Lambda$ to $\oLambda$ in \eqref{eq:cxu1u2}, the strict complementarity condition yields
$\ovlU=\ovlU_1$ and thus $C(\bx)=\mathcal{L}(\bx)$ holds. Therefore, SOSP and {\wSOSP} are identical. This is formally stated in the following proposition.
\begin{proposition}\label{prop:wsosp=sosp}
    Let $\bar{x}$ be a {\rm w}-SOSP of NSDP\,\eqref{eq:main}.
Suppose that the strict complementarity condition holds at $\bx$.
    Then, 
$\bar{x}$ is an SOSP of NSDP\,\eqref{eq:main}.
\end{proposition}
Checking condition \cref{def:SOSPwithoutSC} is computationally intractable. Even for NLPs, it is known to be in general NP-hard from \cite[Theorem~4]{Murty1987}.
On the other hand, the weak second-order necessary conditions can be checked efficiently by computing a basis of $\mathcal{L}(\bar{x})$ if $\Ker(X(\bx))\neq \{0\}$.
\par
We will show that a sequence generated by the proposed IPM accumulates at {\wSOSP}s, which are SOSPs under the strict complementarity condition by Proposition\,\ref{prop:wsosp=sosp}.
As well, in many articles on constrained NLPs,
convergence to SOSPs is established through
    {\wSOSP} and the strict complementarity condition.
For example, see \cite{Auslender1979,FacchineiLucidi1998,Conn2000,Hinder2020}.

\section{Proposed method}\label{chapter:algo}
In this section, we present an IPM, named \hyperref[alg:outer]{Decreasing-NC-PDIPM},
with convergence to SOSPs of NSDP\,\eqref{eq:main}, where each iteration point $(x,Z)\in \euclid^n\times \symmetric^m$ satisfies $X(x)\in \symmetric^m_{++}$ and $Z\in \symmetric^m_{++}$.
\hyperref[alg:outer]{Decreasing-NC-PDIPM} executes an inner algorithm named \hyperref[alg:inner]{Fixed-NC-PDIPM} every iteration.
A core of both the \tteal{algorithms} is the primal-dual merit function $\psi_{\mu,\nu}: \feasi \times \symmetric^m_{++} \rightarrow \euclid$ defined as
\begin{equation}
    \psimunu{x,Z} \coloneqq \FPB(x) + \nu \FBC(x, Z),\label{eq:merit}
\end{equation}
where
\begin{eqnarray*}
  &\FPB(x) \coloneqq f(x) - \mu \logdet{X(x)},&\\
  &\FBC(x, Z) \coloneqq \langle X(x),Z \rangle - \mu \logdet{X(x)} - \mu \logdet{Z},&
\end{eqnarray*}
  and  $\mu, \nu >0$.\footnote{
The subscripts ``PB'' and ``BC'' represent ``primal-barrier'' and ``barrier complementarity'', respectively in abbreviated form.}
This merit function was first presented in~\cite{Yamashita2012} for the NSDP.

\hyperref[alg:outer]{Decreasing-NC-PDIPM} is outlined as follows: 
In order to approach an SOSP of the NSDP, 
\hyperref[alg:outer]{Decreasing-NC-PDIPM} generates an $\varepsilon$-SOSP$(\mu, \nu)$, which will be defined precisely shortly,
while driving $\varepsilon,\mu,\nu >0$ to zeros. 
For each $\varepsilon,\mu,\nu >0$, 
$\varepsilon$-SOSP$(\mu, \nu)$ is an approximate SOSP for $\min_x\psimunu{x,Z}$ whose \tteal{approximation} degree is represented with $\varepsilon>0$.
Each $\varepsilon$-SOSP$(\mu, \nu)$ is computed with \hyperref[alg:inner]{Fixed-NC-PDIPM} taking advantage of a negative-curvature direction, that is, a direction of an eigenvector corresponding to the minimum negative eigenvalue of $\nabla^2_{xx}\psimunu{x, Z}$.
\par
Before defining $\varepsilon$-SOSP$(\mu, \nu)$, we give some formulae related to derivatives of $\psimunu{x, Z}$. First, the partial derivatives of $\psi_{\mu,\nu}$ with respect to $x$ and $Z$ are given as
\begin{align}
    \nabla_x \psimunu{x, Z}
                              & = \nabla f(x)
    - \Aast\paren{x}\paren{\mupn \invX{x}- \nu Z},\label{eq:gradfx}               \\
    \nabla_{Z} \psimunu{x, Z} & = \nu \paren{X(x) - \mu Z^{-1}}.\label{eq:gradfZ}
\end{align}
By differentiating~\cref{eq:gradfx} further, we have for each $i,j$,
\begin{align}
    \paren{\nabla^2_{xx} \psimunu{x, Z}}_{ij}
     & =(\nabla^2_{xx} f(x))_{ij}-\trace{\paren{\mupn \invX{x}- \nu Z}
        \pdv{X}{x_i}{x_j}\paren{x}}                                    \nonumber \\
     & \quad +\mupn
    \trace{
        \Acali\paren{x} \invX{x} \Acalj\paren{x}  \invX{x}}.
    \label{eq:hessofpsi}
\end{align}
Let
\begin{equation}
    \Lambda \coloneqq  (1+\nu)\mu \invX{x} - \nu Z \label{eq:Lambda}.
\end{equation}
We obtain $\nabla_{x} L(x, \Lambda)\allowbreak = \allowbreak \nabla_x \psimunu{x, Z}$ and
\begin{equation}
    (\nabla^2_{xx} \psimunu{x, Z})_{ij} = (\nabla^2_{xx} L(x, \Lambda))_{ij}  +\mupn
    \trace{\Acali\paren{x} \invX{x} \Acalj\paren{x}  \invX{x}}
    \label{eq:HessofLandpsimu}
\end{equation}
by combining~\cref{eq:Lambda} with the derivatives of $\psi_{\mu,\nu}$. The above equations suggest that $\psi_{\mu,\nu}$ is a surrogate for the Lagrangian~$L$. In~\cref{chapter:convergence_outer}, we will reveal the relationship between $\psi_{\mu,\nu}$ and $L$ more precisely.
\par
Now, we define an $\varepsilon$-SOSP$(\mu,\nu)$.
\begin{definition}[$\varepsilon$-SOSP($\mu,\nu$)]
    Given parameters $\varepsilon\coloneqq(\varepsilon_g, \varepsilon_\mu, \varepsilon_H), \mu, \nu>0$, we call $(\xbar,\Zbar) \in \mathcal{X}\times \symmetric^m_{++}$ an $\varepsilon$-SOSP$(\mu, \nu)$ if it holds that
    \begin{subequations}\label{eq:aprox}
        \begin{align}
            \norm{\nabla_x \psimunu{\xbar, \Zbar}}  & \leq   \varepsilon_g\paren{1+\mu\fnorm{X(\xbar)^{-1} }+\fnorm{\Zbar}},\label{eq:apsta} \\
            \fnorm{\nabla_Z \psimunu{\xbar, \Zbar}} & \leq \varepsilon_\mu (1+\mu \fnorm{{\Zbar}^{-1}}),
            \label{eq:apcomple}  \\
            \mineig{\hesspsimu{\xbar, \Zbar}}       & \geq -\varepsilon_H \paren{1+\mu \fnorm{X(\xbar)^{-1}} + \fnorm{\Zbar}}^2.
            \label{eq:apsecond}
        \end{align}
    \end{subequations}
\end{definition}
In fact, by driving $(\varepsilon_g, \varepsilon_\mu,\varepsilon_H)$ and $(\mu, \nu)$ to zeros, it will be shown that
a sequence of $\varepsilon$-SOSP($\mu, \nu$) accumulates at SOSPs of NSDP\,\eqref{eq:main} (cf. \cref{theorem:SOSP}).
Thus, an $\varepsilon$-SOSP($\mu, \nu$) can be regarded as an approximation to an SOSP of NSDP\,\eqref{eq:main}.
\subsection{Decreasing-NC-PDIPM}\label{section:algo_outer}
\hyperref[alg:outer]{Decreasing-NC-PDIPM} produces a sequence of $\varepsilon$-SOSP$(\mu, \nu)$ with decreasing
$\varepsilon,\mu$, and $\nu$ to zeros.
Each $\varepsilon$-SOSP$(\mu, \nu)$ is computed by means of~\hyperref[alg:inner]{Fixed-NC-PDIPM}.
We denote each iteration by $k$ and define
$$
\Xk \coloneqq X(\xk).
$$
  We take continuous functions $\heps_g(\mu), \heps_\mu(\mu), \heps_H\paren{\mu}$, and $\hnu(\mu)$ such that, for each $\mu>0$, the function values are positive and, as $\mu \rightarrow 0$,
  \begin{equation}
  \heps_g(\mu) \rightarrow 0,\ \heps_\mu(\mu) \rightarrow 0,\ \heps_H(\mu) \rightarrow 0,\ \hnu(\mu) \rightarrow 0.\label{eq:parameter}
\end{equation}
The formal description of~\hyperref[alg:outer]{Decreasing-NC-PDIPM} is given in~\cref{alg:outer}.
\begin{remark}
    The proposed algorithm gradually decreases the weight parameter $\nu$ in $\pmunu$ to zero, whereas it is fixed in many PDIPMs\,\cite{Yamashita2012, Okuno2020, kato2015interior, Yamakawa2014}.
    This device is required for technical reasons concerning the proofs. For example, see Proof of~\eqref{eq:FJdualf} in~\Cref{theorem:FJ} and \eqref{eq:Nbounded} in~\cref{theorem:SOSP}.
\end{remark}
\begin{algorithm}[H]
    \DontPrintSemicolon
    \caption{Decreasing Negative Curvature Primal-Dual Interior Point Method (Decreasing-NC-PDIPM)}
    \label{alg:outer}
    \KwIn{
      $(\xinit, Z_1)\in \feasi \times  \symmetric^m_{++}$, $\heps_g(\cdot), \heps_\mu(\cdot), \heps_H(\cdot), \hnu(\cdot)$,
       $\mu_1, L_0, L_1, L_2>0$}
    $k\leftarrow 1$\;
    \While{not converge}{
        Set $\mu_{k+1}>0$ so that $\mu_{k+1}<\mu_k$ and $\lim_{k\to\infty}\mu_k=0$} \;
        $(\mu,\nu, \varepsilon_g, \varepsilon_\mu,\varepsilon_H)\leftarrow (\mukp,\hnu(\mukp), \heps_g(\mukp), \heps_\mu(\mukp),\heps_H(\mukp))$\;
        $(x^{k+1}, Z_{k+1}) \leftarrow$  Output of Fixed-NC-PDIPM($\mu,\nu, \varepsilon_g, \varepsilon_\mu,\varepsilon_H$)
            started from $\paren{x^{k}, Z_{k}}$\;
        $k \leftarrow k+1$
\end{algorithm}
\subsection{Fixed-NC-PDIPM}\label{section:algo_inner}
We describe \hyperref[alg:inner]{Fixed-NC-PDIPM} formally in~\cref{alg:inner}, where we define
\begin{equation}
  \Xell \coloneqq X(\xell)
  \label{eq:xell}
\end{equation}
for each $\ell$.
Given the parameters $(\varepsilon,\mu, \nu)$ and an initial point from \hyperref[alg:outer]{Decreasing-NC-PDIPM},
\hyperref[alg:inner]{Fixed-NC-PDIPM} computes an $\varepsilon$-SOSP$(\mu, \nu)$ by searching $\feasi\times \mathbb{S}^m_{++}$
\par
Each of the Updates in \hyperref[alg:inner]{Fixed-NC-PDIPM} employs a specific direction:
\begin{itemize}[label=\textbullet]
    \item $-\frac{\iotamin}{\iotamax^2} \calHell{\nabla_{Z} \psi_{\mu,\nu}}$ in \autoref{alg:dual},
    \item $-\frac{\hmin}{\hmax^2}\Hell \nabla_{x} \psi_{\mu,\nu}$ in \autoref{alg:primal},
    \item an eigenvector corresponding to $\mineig{\nabla^2_{xx} \psi_{\mu,\nu}}$ in \autoref{alg:negative_curvature},
\end{itemize}
where $\iotamin,\iotamax,\hmin,\hmax>0$ are prefixed constants,
$\mathcal{H}_{\ell}$ is a {\it symmetric} linear operator from $\symmetric^m$ to $\symmetric^m$,
and $H_{\ell}$ is an $n \times n$ positive definite symmetric matrix satisfying
\begin{align}
     & 0 < \iotamin \leq \langle D, \allowbreak \calHell{D} \rangle \allowbreak \leq \iotamax,\ \forall D \in \symmetric^m \mbox{ with }\|D\|_{\rm F}=1, \label{eq:calHell}\\
     & O \prec \hmin I_n \preceq \Hell \preceq \hmax I_n.
     \label{eq:Hell}
\end{align}
In order to obtain a faster convergence, it is desirable to set
the inverse of the Hessian or their approximation as $\mathcal{H}_{\ell}$ and $H_{\ell}$.
When  $\iotamin=\iotamax=\hmin=\hmax=1$, $H_{\ell}=I_m$ and $\calHell\cdot$ is the identity mapping, 
the directions in Updates~1 and~2
reduce to the usual partial derivatives.
\par
The three directions are devised so as to ensure conditions~\eqref{eq:apsta},~\eqref{eq:apcomple}, and~\eqref{eq:apsecond}, respectively.
If any of those conditions is violated, proceeding along the corresponding direction with the step size $\alpha_{\ell}$, we can
\begin{itemize}
    \item keep the strict feasibility (\cref{proposition:forinner}),
    \item decrease $\psi_{\mu,\nu}$ by more than a certain constant (\cref{lemmabreak:descent}).
\end{itemize}
The step size $\alpha_{\ell}$ is determined according to the \Lipschitz constants $L_0,L_1,$ and  $L_2$. Even when we do not know these constants {\it a priori}, backtracking line search enables us to \tteal{invoke} Updates~1,2, and 3.
In numerical experiments, we
implement this technique.
Actually, the same convergence guarantees hold even if the order of Procedures 1, 2, and 3 is changed in \hyperref[alg:inner]{Fixed-NC-PDIPM}.
However, taking into consideration the computation cost of eigenvalue decomposition, it would be favorable to prioritize the use of the gradient directions over the minimum eigenvalue direction.
\setcounter{algocf}{1}
\begin{algorithm}
    \DontPrintSemicolon
    \caption{Fixed Negative Curvature Primal-Dual Interior Point Method with $\mu$, $\nu$, $\varepsilon_g$,
        $\varepsilon_\mu$, $\varepsilon_H$ (Fixed-NC-PDIPM ($\mu,\nu, \varepsilon_g, \varepsilon_\mu,\varepsilon_H$))}
    \label{alg:inner}
    \KwIn{$\xinit, Z_{1}, \mu, \nu, \varepsilon_g, \varepsilon_\mu, \varepsilon_H, L_0, L_1, L_2, \hmin, \hmax, \iotamin, \iotamax$}\tcp{$(\xinit, Z_1)\in \feasi \times \mathbb{S}^m_{++}$}
    \For{$\ell = 1, \ldots$ }{
        \uIf{$\fnorm{\nabla_Z \psimunu{\xell, \Zell}} > \varepsilon_\mu (1+\mu \fnorm{\Zell^{-1}})$\tcp{Procedure 1}}{
            Set $\Zellp$ by \autoref{alg:dual}\;
            $\xellp \leftarrow \xell$
        }
        \uElseIf{$\norm{\nabla_x \psimunu{\xell, \Zell}}>   \varepsilon_g\paren{1+\mu\fnorm{\invXell}+\fnorm{\Zell}}$\tcp{Procedure 2}}{
            Set $\xellp$ by \autoref{alg:primal}\;
            $\Zellp \leftarrow \Zell$
        }
        \uElseIf{$\mineig{\hesspsimu{\xell, \Zell}} < -\varepsilon_H \paren{1+\mu \fnorm{\invXell} + \fnorm{\Zell}}^2$\;\tcp{Procedure 3}}{
            Set $\xellp$ by \autoref{alg:negative_curvature}\;
            $\Zellp \leftarrow \Zell$
        }
        \Else{
            Output $(\xell, \Zell)$
        }
    }
\end{algorithm}

\setcounter{algocf}{0}
\begin{procedure}
    \DontPrintSemicolon
    \caption{Use of the scaled partial derivative w.r.t. $Z$ at () $\paren{\xell,\Zell}$}
    \label{alg:dual}
    \KwIn{$\xell, \Zell, \calHell{\cdot}, \iotamin, \iotamax$}
    \KwOut{$\Zellp$}
    ${d}_{\Zell} \leftarrow -\frac{\iotamin}{\iotamax^2} \calHell{\nabla_{Z} \psi_{\mu,\nu}}$\;
    $l_Z(\Zell) \leftarrow 2\mu\nu \fnorm{\Zell^{-1}}^2 $
    \tcp{``Local'' \Lipschitz constant of $\nabla_Z \psi_{\mu,\nu}$}
    $\innerstepsize \leftarrow \min\left(\frac{1}{l_Z(\Zell)}, \frac{\mineig{\Zell}}{2 \fnorm{{d}_{\Zell}}}\right)$ \tcp{Step size}
    $\Zellp \leftarrow \Zell + \innerstepsize {d}_{\Zell} $
\end{procedure}

\begin{procedure}
    \DontPrintSemicolon
    \caption{Use of the scaled partial derivative w.r.t. $x$ at () $\paren{\xell, \Zell}$}
    \label{alg:primal}
    \KwIn{$\xell, \Zell, H_{\ell}, L_0, L_1, \hmin, \hmax$}
    \KwOut{$\xellp$}
    ${d}_{\xell} \leftarrow -\frac{\hmin}{\hmax^2}\Hell \nabla_{x} \psi_{\mu,\nu}$\;
    $l_x(\xell, \Zell) \leftarrow L_1 + \nu L_1\fnorm{\Zell}+ 2\mupn L_0^2\fnorm{\invXell}^2+ \mupn L_1 \fnorm{\invXell}$
    \tcp{``Local'' \Lipschitz constant of $\nabla_x \psi_{\mu,\nu}$}
    $\innerstepsize \leftarrow
        \min\left(\frac{\mineig{\Xell}}
        {2 L_0 \nnorm{\innerdirection}},\frac{1}{l_x(\xell, \Zell)}\right)$ \tcp{Step size}
    $\xellp \leftarrow \xell + \innerstepsize \innerdirection$
\end{procedure}
\begin{procedure}
    \DontPrintSemicolon
    \caption{Use of the negative curvature w.r.t. $x$ at () $\paren{\xell, \Zell}$}
    \label{alg:negative_curvature}
    \KwIn{$\xell, \Zell, L_0, L_1, L_2$}
    \KwOut{$\xellp$}
    $\innerdirection \leftarrow$ The normalized eigenvector corresponding to $\lambda_{\min}\paren{\nabla_{xx}^2\psimunu{\xell,\Zell}}$\;
    \tcp{``Local'' \Lipschitz constant of $\nabla^2_{xx} \psi_{\mu,\nu}$}
    $l_{xx}(\xell, \Zell) \leftarrow L_2+\nu L_2 \fnorm{\Zell}+
        \mupn\paren{ L_2 \fnorm{\invXell} +
            4 L_1 L_0\fnorm{\invXell}^2 + 6 L_0^3 \fnorm{\invXell}^3}$\;
   \tcp{Reverse the direction in order to ensure $\innerdirection$ is a descent direction}
    \If{$\innerdirection^\top \nabla_x \psimunu{\xell, \Zell} > 0$}{
        $\innerdirection \leftarrow - \innerdirection$}
    $\innerstepsize \leftarrow \min\left(-\frac{2 \mineig{\hesspsimu{\xell, \Zell}}}{l_{xx}(\xell, \Zell)},\frac{\mineig{\Xell}}{2L_0 \norm{\innerdirection}}\right)$
    \tcp{Step size}
    $\xellp \leftarrow \xell + \innerstepsize \innerdirection$
\end{procedure}

\section{Convergence to SOSPs}
\label{chapter:convergence_outer}
In this section, supposing that \hyperref[alg:outer]{Decreasing-NC-PDIPM} produces an infinite sequence of $\varepsilon$-SOSP$(\mu,\nu)$,
we show that the sequence accumulates at SOSPs in the following manner.
First, \cref{section:FJ} asserts the convergence to FJ points under certain reasonable assumptions.
Second, in~\cref{section:KKT}, the additional assumption of the MFCQ yields the convergence to KKT points. Finally, in \cref{section:SOSP}, we prove that these KKT points are equivalent to {\wSOSP}s.
\par
Since each iterate of~\hyperref[alg:outer]{Decreasing-NC-PDIPM} is 
an $\varepsilon_k$-SOSP($\mu_k, \nu_k$), the following equations hold:
\begin{subequations}
    \begin{align}
        \Xk                                 & \in \symmetric^m_{++},\label{eq:innterfeasi}                                            \\
        \Zk                                 & \in \symmetric^m_{++},                                                                  \\
        \norm{\nabla_x \psimunu{\xk, \Zk}}  & \leq   \heps_g(\mu_k)\paren{1+\mu_k\fnorm{\invXk}+\fnorm{\Zk}}, \label{eq:outsta} \\
        \fnorm{\nabla_Z \psimunu{\xk, \Zk}} & \leq \heps_{\mu}(\mu_k) (1+\mu_k \fnorm{\invZk}),
        \label{eq:outcomple}                                                                                                          \\
        \mineig{\hesspsimu{\xk, \Zk}}       & \geq -\heps_H(\mu_k) \paren{1+\mu_k \fnorm{\invXk} + \fnorm{\Zk}}^2.
        \label{eq:outsecond}
    \end{align}
\end{subequations}
\begin{remark}
The subsequent analysis on convergence to SOSPs focuses only on the behavior of $\varepsilon_k$-SOSP($\mu_k, \nu_k$) as $k$ tends to $\infty$. The analysis never requires \cref{assumption:Lipschtiz}, which is the assumption of the \Lipschitz constants of the functions. 
\end{remark}
\subsection{Convergence to FJ points}
\label{section:FJ}
In this section, we will prove that any accumulation points of~\hyperref[alg:outer]{Decreasing-NC-PDIPM} are FJ points under the following assumption.
\begin{assumption}
    \leavevmode
    \begin{enuminasm}
        \item The sequence $\sbra{\xk}$ is bounded. \label{assumption:bounded}
        \item
        The sequence $\sbra{\Zk - \muk \invXk}$ is bounded.
        \label{assumption:zmuinvxbounded}
    \end{enuminasm}
    \label{assumption:forFJ}
\end{assumption}
The first assumption is often found in many articles for NLPs or NSDPs, for example,~\cite{Moguerza2003,Yamashita2012,Yamashita2020}.  
In contrast, the second one is actually closely related to \cref{eq:outcomple}, and may be peculiar to the analysis in this paper.   
On the central-path-like set formed by $0$-SOSP$(\mu,\nu)$, 
$Z = \mu X^{-1}$ holds from \cref{eq:outcomple}. In view of this fact, the assumption is interpreted as 
that the sequence is not too distant from the path-like set.  
\begin{theorem}[convergence to FJ points]
    \label{theorem:FJ}
    Suppose that~\cref{assumption:forFJ} holds.
    Then,
    $\sbra{(\xk, \lambda_k, \Omega_k)}$ with
    \begin{align*}
        \lambda_k \coloneqq  \frac{1}{1+\muk \fnorm{\invXk}+\fnorm{\Zk}},\;
        \Omega_k \coloneqq \frac{\Lambda_k}{1+\muk \fnorm{\invXk}+\fnorm{\Zk}}
    \end{align*}
    is bounded and any accumulation points satisfy the FJ conditions, where $\Lambda_k=  \mukpn \invXk- \hnu(\muk) \Zk$ as defined in~\cref{eq:Lambda}.
\end{theorem}
\begin{proof}
    Note that from \cref{eq:parameter,assumption:zmuinvxbounded} and $\lim_{k\to\infty}\mu_k=0$,
        \begin{equation}
            \lim_{k\to\infty}\hnu(\muk) \fnorm{\muk \invXk - \Zk}=0
            \label{prop:conv0}
        \end{equation}
    holds.
    Since $\{\xk\}$ is bounded from~\cref{assumption:bounded} and $X$ is continuous, $\{\Xk\}$ is bounded.
    Next, we prove the boundedness of $\{(\Omega_k,\lambda_k)\}$. From the definition of $\Omega_k$,  we have
    \begin{align*}
        \fnorm{\Omega_k}
         & =\frac{\fnorm{\muk \invXk + \hnu(\muk) \paren{\muk \invXk - \Zk}}}{1+ \muk \fnorm{\invXk} +\fnorm{\Zk}} \\
         & \leq
        \frac{\muk \fnorm{\invXk}}{1+ \muk \fnorm{\invXk} +\fnorm{\Zk}}
        + \frac{\hnu(\muk) \fnorm{\muk \invXk - \Zk}}{1+ \muk \fnorm{\invXk} +\fnorm{\Zk}}                         \\
         & \leq 1
        + {\hnu(\muk) \fnorm{\muk \invXk - \Zk}}.
    \end{align*}
    By taking $\limsup$ of $\fnorm{\Omega_k}$ and the rightmost hand and using \eqref{prop:conv0}, we see that
    $\limsup_{k\rightarrow \infty} \allowbreak \fnorm{\Omega_k}\le 1$, implying the boundedness of
    $\{ \Omega_k \}$. Since $ 0< \lambda_k = \frac{1}{1+\muk \fnorm{\invXk}+\fnorm{\Zk}} \allowbreak \leq 1,$
    $\sbra{\lambda_k}$ is bounded.
    Consequently, we obtain the boundedness of $\{(x^k,\allowbreak \lambda_k,\allowbreak \Omega_k,\allowbreak X_k)\}$.
    \par
    Let an arbitrary accumulation point of $\bigl\{\bigl(\xk, \Xk, \lambda_k, \Omega_k,\bigr)\bigr\}$ be $\bigl(\xbar,\allowbreak\overline{X}, \allowbreak\overline{\lambda},\allowbreak \overline{\Omega}\bigr).$ To prove the convergence to FJ points, let us check conditions~\cref{eq:FJfeasi,eq:FJdualf,eq:FJcomp,eq:FJnonzero,eq:FJstationary} one by one.
    \begin{proof}[Proof of~\cref{eq:FJstationary}]
        We have,  from condition~\cref{eq:outsta},$
            \norm{\nabla_x \psimunu{\xk, \Zk}} \leq   \heps_g(\mu_k)\allowbreak \bigl(1\allowbreak+\allowbreak \mu_k \allowbreak \bigl\|\invXk\bigr\|_{\mathrm{F}}\allowbreak+\allowbreak\fnorm{\Zk}\bigr).$
        By substituting~\cref{eq:gradfx} into this inequality and dividing both sides with
        $1+\mu_k\bigl\|\invXk\bigr\|_{\mathrm{F}}\allowbreak+\allowbreak\fnorm{\Zk}$,
        we obtain \[
            \|\lambda_k \nabla f(\xk) - \Aast\paren{\xk} \Omega_k\|
            \leq \heps_g(\muk)\xrightarrow[k \rightarrow \infty]{} 0\] and thus~\cref{eq:FJstationary} is ensured at $\bigl(\xbar, \allowbreak\overline{\lambda},\allowbreak \overline{\Omega}\bigr).$
    \end{proof}
    \begin{proof}[Proof of~\cref{eq:FJcomp}]
        From $\lim_{k\to\infty}\mu_k=0$, \cref{eq:parameter}, and \cref{assumption:zmuinvxbounded}, we have
        \begin{align*}
            \fnorm{\Lambda_k \Xk} & =\fnorm{\muk I_m
                + \hnu(\muk) \paren{\muk I_m-  \Zk \Xk}}   \\
                                  & \leq \fnorm{\muk I_m}
            + \hnu(\muk) \fnorm{\Xk}\fnorm{\muk \invXk -\Zk } \xrightarrow[k \rightarrow \infty]{} 0.
        \end{align*}
        Therefore we have
        \begin{equation*}
            \fnorm{\Xk\Omega_k} =
            \frac{\fnorm{\Xk\Lambda_k}}{1+\muk \fnorm{\invXk}+\fnorm{\Zk}}\leq \fnorm{\Xk\Lambda_k} \xrightarrow[k \rightarrow \infty]{} 0
            \nonumber
        \end{equation*}
        and thus~\cref{eq:FJcomp} holds at $\bigl(\overline{X}, \overline{\Omega}\bigr).$

    \end{proof}
    \begin{proof}[Proof of~\eqref{eq:FJfeasi}]
        From condition~\eqref{eq:innterfeasi}, any accumulation points of $\{\Xk\}$ satisfy~\eqref{eq:FJfeasi}.
    \end{proof}
    \begin{proof}[Proof of~\eqref{eq:FJdualf}]
        Since $\lambda_k\ge 0$ at each $k$, its accumulation point is also nonnegative.
            We next show that any accumulation point of $\{\Omega_k\}$ is positive semidefinite. First, it follows that
        \begin{align*}
            \frac{\hnu(\muk) \mineig{
                    \muk \invXk - \Zk
                }}{1+ \muk \fnorm{\invXk} +\fnorm{\Zk}}
             & \leq \hnu(\muk) \fnorm{\muk \invXk - \Zk}
        \end{align*}
        and thus
        the left-hand side converges to zero as $k\to\infty$ from
            \cref{prop:conv0}.
        Since $\mineig{A} + \mineig{B} \leq \mineig{A+B}$ for $A,B \in \symmetric^m_{+}$ (see \cite[Theorem 4.3.1]{Horn2012}), we have, from the above inequality,
        \begin{align*}
             & \liminf_{k\rightarrow \infty}  \mineig{\Omega_k}                                                                                              \\
             & \geq  \liminf_{k\rightarrow \infty} \frac{ \mineig{\muk \invXk}}{1+ \muk \fnorm{\invXk} +\fnorm{\Zk}} +  \liminf_{k\rightarrow \infty}\left(
            \hnu(\muk)  \frac{\mineig{
                    \muk \invXk - \Zk
                }}{1+ \muk \fnorm{\invXk} +\fnorm{\Zk}}                                                                                              \right) \\
             & \geq 0.
        \end{align*}
        Therefore~\cref{eq:FJdualf} is ensured at $\overline{\Omega}$.
    \end{proof}
    \begin{proof}[Proof of~\cref{eq:FJnonzero}]
        \begin{align*}
            \lambda_k + \fnorm{\Omega_k}
             & = \frac{1+\fnorm{\muk \invXk+\hnu(\muk) \paren{\muk \invXk - \Zk}}}{1+\muk \fnorm{\invXk}+\fnorm{\Zk}}
            \\
             &
            \geq \frac{1+\fnorm{\muk \invXk}}{1+\muk \fnorm{\invXk}+\fnorm{\Zk}}
            -\frac{\hnu(\muk) \fnorm{\muk \invXk - \Zk}}{1+\muk \fnorm{\invXk}+\fnorm{\Zk}}
        \end{align*}
        To obtain a contradiction, assume that there exists a subsequence $\bigl\{\lambda_k \allowbreak+ \allowbreak\fnorm{\Omega_k}\bigr\}_{k\in \mathcal{S}}$ converging to zero. Taking the limit of the above inequality and using \eqref{prop:conv0} again yield
        \[
        0=
            \lim_{k \in \mathcal{S} \rightarrow \infty} \left(\lambda_k + \fnorm{\Omega_k}\right)
            \geq  \lim_{k \in \mathcal{S}\rightarrow \infty}
            \frac{1+\fnorm{\muk \invXk}}{1+\muk \fnorm{\invXk}+\fnorm{\Zk}},\\
       \]
        which entails 
        \begin{equation}
            \lim_{k\in \mathcal{S} \rightarrow \infty}
            \frac{1+\fnorm{\muk \invXk}}{1+\muk \fnorm{\invXk}+\fnorm{\Zk}} = 0.
            \label{eq:limzero}
        \end{equation}
        We consider the two cases where $\{Z_k\}_{k\in \mathcal{S}}$ is bounded and unbounded and then derive contradictions in both two cases.
        \par
        If $\{Z_k\}_{k\in \mathcal{S}}$ is bounded, by \cref{assumption:zmuinvxbounded}, $\{\mu_k X_k^{-1}\}_{k\in \mathcal{S}}$ is also bounded.
        Then, the left-hand side in \eqref{eq:limzero} must be positive because
        it is bounded from below by $1/\limsup_{k\in \mathcal{S}\to\infty}\bigl(1\allowbreak+\allowbreak\muk \fnorm{\invXk}\allowbreak+\allowbreak \fnorm{\Zk}\bigr)$, but this contradicts \eqref{eq:limzero}.
        \par
        Next, consider the case where $\{Z_k\}_{k\in \mathcal{S}}$ is unbounded. Notice that it holds that
        \begin{align}
            \frac{1+\fnorm{\muk \invXk}}{1+\muk \fnorm{\invXk}+\fnorm{\Zk}}
            \ge \frac{1+\left| \fnorm{\muk \invXk-Z_k} - \fnorm{Z_k}\right|}{1+\fnorm{\mu_k\invXk-\Zk}+2\fnorm{\Zk}}\eqqcolon v_k, \label{al:2124}
        \end{align}
        where the inequality follows by 
        the fact $\fnorm{A+B}^2\ge (\fnorm{A}-\fnorm{B})^2$ for any $A,B\in \symmetric^m$ and the substitution $(A,B)=(\mu_k\invXk-\Zk,\Zk)$. 
        \Cref{assumption:zmuinvxbounded} and the unboundedness of $\sbra{\|Z_k\|_F}_{k\in \mathcal{S}}$ yield $\limsup_{k \in \mathcal{S} \to\infty}v_k\ge 1/2$. Meanwhile, \eqref{eq:limzero}, \eqref{al:2124}, and $v_k \geq 0$ lead to
        $\lim_{k\to\infty}v_k=0$, giving a contradiction again.
        So, \Cref{eq:FJnonzero} holds at $\bigl(\overline{\lambda}, \overline{\Omega}\bigr)$.
    \end{proof}
    The proof is complete.\qed
\end{proof}
\subsection{Convergence to KKT points}
\label{section:KKT}
In this section, we prove that any accumulation points of a sequence generated by~\hyperref[alg:outer]{Decreasing-NC-PDIPM} are KKT points by assuming the MFCQ in addition to~\cref{assumption:forFJ}.
\begin{assumption}[MFCQ]
    \label{assumption:MFCQ}
    The MFCQ holds at any accumulation point of $\{\xk\}$.
\end{assumption}
\begin{theorem}[convergence to KKT points]
    Suppose that \Cref{assumption:forFJ,assumption:MFCQ} hold. Then
    $\{(\xk, \allowbreak \Lambda_k)\}$ is bounded and any accumulation points of $\{(\xk, \allowbreak \Lambda_k)\}$ satisfy the KKT conditions. In addition, $\{\Zk\}$ and $\{\muk \invXk\}$ are bounded.
    \label{theorem:KKT}
\end{theorem}
\begin{proof}
    Let $\{\lambda_k\}$ and $\{\Omega_k\}$ be the same ones as defined in~\cref{theorem:FJ}.
    From \cref{theorem:FJ}, $\{\paren{\xk, \lambda_k, \Omega_k}\}$ is bounded and furthermore, by letting
    $(\bx,\blambda,\overline{\Omega})$ be its arbitrary accumulation point, $(\bx,\blambda,\overline{\Omega})$ satisfies the FJ conditions.
    Under the MFCQ, $\bx$ is a KKT point and $\blambda\neq 0$.
    This fact together with
    the boundedness of $\{\Omega_k\}$ yields the boundedness of $\{\Lambda_k\}=\{ \lambda_k^{-1}\Omega_k\}$.
    Therefore, any accumulation point of $\{(x^k,\Lambda_k)\}$ satisfies the KKT conditions. We thus obtain the former assertion.
    \par
    Since an arbitrary accumulation point of $\{\lambda_k\}$ is nonzero,
    we see that $\sbra{\Zk}$ and $\sbra{\muk \invXk}$ are bounded by recalling the definition of $\lambda_k$. The proof is complete.
    \qed
\end{proof}
\begin{assumption}
    The sequence $\sbra{\frac{\heps_\mu (\muk)}{\hnu(\muk) \muk}}$ converges to zero.
    \label{assumption:compleconv}
\end{assumption}
\begin{theorem}\label{theorem:Zlambda}
    Suppose that~\Cref{assumption:MFCQ,assumption:compleconv} hold. Then we have $\allowbreak
        \bigl\|\muk\allowbreak \invXk - \Zk\bigr\|_{\mathrm{F}} \xrightarrow[k\rightarrow \infty]{} 0,\allowbreak  \fnorm{\Lambda_k - \Zk}\xrightarrow[k\rightarrow \infty]{} 0.$
\end{theorem}
\begin{proof}
    From~\eqref{eq:outcomple} and~\eqref{eq:gradfZ}, we have
    \[
        \fnorm{\Xk - \muk \Zk^{-1}} \leq \frac{\heps_{\mu}(\muk)}{\hnu(\muk)} (1+\muk \fnorm{\Zk^{-1}}),
    \]
    and thus get
    \begin{align*}
        \fnorm{\muk \invXk - \Zk} & =
        \fnorm{\invXk\paren{\muk \Zk^{-1} - \Xk}\Zk}                                                                                                                                      \\
                                  & \leq \fnorm{\invXk} \fnorm{\muk \Zk^{-1} - \Xk}\fnorm{\Zk}                                                                                            \\
                                  & \leq \frac{\heps_\mu(\muk)}{\hnu(\muk)\muk} \fnorm{\muk \invXk} \fnorm{\Zk}\paren{1+\muk \fnorm{\Zk^{-1}}}\xrightarrow[k\rightarrow \infty]{} 0,
    \end{align*}
    where the last convergence follows from~\Cref{assumption:compleconv} and the boundedness of $\sbra{\muk \invXk}$ and $\sbra{\Zk}$ from \cref{theorem:KKT}. We also have $
        \fnorm{\Lambda_k - \Zk}
        \allowbreak=\allowbreak(1+\hnu(\muk))\fnorm{\muk \invXk - \Zk }\allowbreak \xrightarrow[k\rightarrow \infty]{}0.$
    The proof is complete.
    \qed
\end{proof}
\subsection{Convergence to SOSPs}
\label{section:SOSP}
Before presenting the theorem regarding the convergence to {\wSOSP}s, we define some terminology.
Let $\mathcal{D}^p$ and $\mathcal{D}^p_{++(+)}$ be
the set of $p\times p$ diagonal matrices and  the set of  $p\times p$ positive (semi) definite diagonal matrices, respectively.
    We also define
    $\mathcal{O}^{p\times q}\coloneqq\{U\in \euclid^{p\times q}\mid U^{\top}U =I_q\}$. If $p=q$, we abbreviate $\calO^{p\times p}$ as $\calO^{p}$.
    Under~\Cref{assumption:bounded},~$\sbra{\Xk}$ and $\{x^k\}$ are bounded and thus have accumulation points, say $X_{\ast}$ and $\xast$, respectively. There exists a subsequence $\calKo$ such that
    \[\lim_{k\in\calKo\to\infty}x^k=\xast,\ \lim_{k \in \calKo \rightarrow \infty} \Xk=\Xast.\]
In what follows, we introduce some specific sequences related to the subsequence $\calKo$.
For each $k\in\calKo$, an eigenvalue decomposition of $\Xk \in \symmetric^m_{++}$ is expressed as
\begin{equation*}
    \Xk =  \Vk \Dk {\Vk}^\top,\ \Vk\in \calO^m,\ \Dk\in \calD^m_{++}.
\end{equation*}
Since $\sbra{\Vk}$ and $\sbra{\Dk}$ are bounded, they have accumulation points, say $\Vast\in \calO^m$ and $\Dast\in \calD^m_{+}$.
We have $\Xast= \Vast \Dast \Vast^\top.$
Let \[\rast \coloneqq \rank (\Xast)\] and $\Dastp$ be the unique submatrix of $\Dast$
such that it is in $\mathcal{D}^{\rast}_{++}$.
Extracting appropriate $\rast$ columns of $\Vast$, written as $\Vastp\in\calO^{m \times \rast}$,  we can rewrite $\Xast$ as $X_\ast = \Vastp \Dastp \transpose{\Vastp}$.
Moreover, let $\Vastz \in \calO^{m\times (m-\rast)}$ be the remaining columns of $\Vast$ after excluding $\Vastp$ from $\Vast$.
Without loss of generality, we assume
\begin{equation}
    \Dast=\begin{bmatrix}\Dastp&O\\O&O\end{bmatrix},\ \Vast = \left[\Vastp,\Vastz\right],\ \lim_{k\in \calKo\to\infty}(\Vk,\Dk)=(\Vast,\Dast).\label{eq:without}
\end{equation}
Notice that $\Vastz$ spans $\Ker(X_\ast)$, while $\Vastp$ does its orthogonal complement space.
In relation to the above decomposition of $X_\ast$,
there exist some matrices $\Dkp\in\calD^{\rast}_{++}, \Vkp\in \calO^{m\times \rast},\Dkz\in \calD^{m-\rast}_{++}$ and $\Vkz\in\calO^{m\times (m-\rast)}$
such that
\[
    \Xk = \Vkp \Dkp \transpose{\Vkp} + \Vkz \Dkz \transpose{\Vkz}
\] and
\[
    \Dk=\begin{bmatrix}\Dkp&O\\O&\Dkz\end{bmatrix},\ \Vk = [\Vkp,\Vkz].
\]
Notice that
\begin{equation}
    \lim_{k\in\calKo\to\infty}\left(\Dkp,\Dkz,\Vkp,\Vkz\right) = \left(\Dastp,\Dastz,\Vastp,\Vastz\right)
    \label{eq:limitmany}
\end{equation}
by comparing the above equations and \eqref{eq:without}.
Furthermore, we define 
$$\Xkp \coloneqq \Vkp \Dkp \transpose{\Vkp},\ \Xkn \coloneqq \Vkz \Dkz \transpose{\Vkz},$$ leading to
\begin{align}
    \invXkP & =\Vkp \invDkp \transpose{\Vkp},\nonumber            \\
    \invXkN & = \Vkz \invDkz \transpose{\Vkz}, \label{def:invXkN} \\
    \invXk  & = \invXkP + \invXkN. \label{eq:invDecomp}
\end{align}
\par
If $\Ker(\Xast)=\{0\}$, i.e., $x^{\ast}\in \feasi$, the sequences $\{\Vkz\}_{k\in\calKo}$, $\{\Dkz\}_{k\in\calKo}$,
$\{\Xkn\}_{k\in\calKo}$, and the points $\Vastz$ and $\Dastz$ are not well-defined.
Even in that case, the above argument makes sense just by neglecting the terms including them.
\subsubsection{Theorem for convergence to SOSP}
We require an additional assumption on the subsequence $\calKo$.
To begin with,
for each $k\in\calKo$, let $G_k\in \euclid^{\binom{m-\rast+1}{2} \times n}$ be the matrix whose each row is expressed as
\[
\left[{(\Vkz)_p}^\top \mathcal{A}_1 (\xk) (\Vkz)_q  , \ldots,
        {(\Vkz)_p}^\top \mathcal{A}_n (\xk) (\Vkz)_q
        \right]\ (1 \leq p \leq q \leq m-\rast),
\]
where $\Acali(\cdot)$ ($i=1,\ldots,n$) are defined in \cref{subsection:notation}
and moreover, let
\begin{align*}
    \Lastk \allowbreak & \coloneqq \allowbreak
    \Ker(G_k)
    =\Bigl\{ u \in \euclid^n \allowbreak \Bigm|\allowbreak \paren{\Vkz}^\top \Dx{\xk}{u}\Vkz \allowbreak= O \Bigr\}.
\end{align*}
Notice that $\{G_k\}_{k\in\calKo}$ has a limit, denoted by $\Gast$, as $\lim_{k\in \calKo\to\infty}(x^k,\Vkz)=(x^{\ast},\Vastz)$.
It follows that $\mathcal{L}(x^{\ast})=\Ker(\Gast)$ from the definitions of $\calL(\xast)$ and $\Gast$.

If $\Ker(\Xast)=\{0\}$, the above symbols are not well-defined and $\mathcal{L}(\xast)=\euclid^n$.
In this case, we define $\Lastk\coloneqq \euclid^n$ for each $k\in\calKo$.
\begin{assumption}
    If $\Ker(X_{\ast})\neq \{0\}$, then $\rank(G_k)=\rank(\Gast)$ for any sufficiently large $k \in \calKo$.
    \label{assumption:nondegeneracy}
\end{assumption}
Associated with this assumption, we give the following proposition that will play an important role to prove the theorem on convergence to SOSPs.
An analogous proposition is found in \cite[Lemma 3.1]{Andreani2007}, and
the following one can be proved in a manner quite similar to this lemma, so we omit the proof here.
\begin{proposition}
    Assume that \cref{assumption:nondegeneracy} holds.
    For arbitrarily chosen $\bar{u} \in \calL(\xast)$ with $\calL$ defined in \cref{eq:linearspaceL},
    there exists a sequence $\{u^k\}_{k\in\calKo}$ converging to $\bar{u}$ such that $u^k \in \Lastk$ for each $k\in\calKo$.
    \label{prop:innersemi}
\end{proposition}
Before moving onto the theorem on SOSPs, we make a remark about a sufficient condition for \cref{assumption:nondegeneracy}.
\begin{remark}
    A sufficient condition for \cref{assumption:nondegeneracy} is the nondegeneracy condition (NC)
    \cite[p.86]{wolkowicz2012handbook} at $\xast$,
    a constraint qualification of NSDP\,\eqref{eq:main}.
    Indeed, suppose $\Ker(\Xast))\neq \{0\}$.
    Since $\lim_{k\in\calKo\to\infty}G_k=\Gast$ and the NC means that $\Gast$ is of full row rank, we obtain
    $\rank(G_k)=\rank(\Gast)$ for all $k\in \calKo$ large enough, meaning \cref{assumption:nondegeneracy}.
\end{remark}
\par
Now, we are ready to present the theorem on convergence to SOSPs.
The key ingredient of the proof is that we show that, for any sequence $\{u^k\}_{k\in\calKo}$ with $\uk \in \Lastk$,
$$\liminf_{k\in\calKo\rightarrow \infty} {\uk}^\top\allowbreak \HessLagrange{\xk, \Lambda_k} \allowbreak \uk \allowbreak+\allowbreak \quadf{\widehat{H}(\xk, \Lambda_k)}{\uk} \geq 0,
$$
where
\begin{equation}
    \widehat{H}(\xk, \Lambda_k)_{ij} \coloneqq 2\trace{\Acali(\xk) \invXkP \Acalj(\xk) \Lambda_k}\label{eq:Hhat}
\end{equation}
for each $i,j$.
In fact, the second term ${\uk}^\top\allowbreak \widehat{H}(\xk, \Lambda_k)\allowbreak \uk$ converges to the sigma-term at an accumulation point of $\{x^k\}$, which leads us to the weak second-order necessary conditions.
\begin{theorem}[convergence to {\wSOSP}s]
    Suppose that~\Cref{assumption:forFJ,assumption:MFCQ,assumption:nondegeneracy} hold. Then any accumulation points of $\sbra{\paren{\xk, \Lambda_k}}$ are {\wSOSP}s of NSDP\,\eqref{eq:main}, which are identical to SOSPs when the strict \comple condition holds there.
    \label{theorem:SOSP}
\end{theorem}
\begin{proof}
    \label{proof:SOSP}
        The sequence
        $\{x^k\}$ is bounded due to \cref{assumption:bounded}, and has an accumulation point, say $x^{\ast}$,
        for which the subsequence $\calKo$ is constructed as at the beginning of this subsection.
        Since $\sbra{\Lambda_k}_{k\in\calKo}$ is bounded by \cref{theorem:KKT}, we have an accumulation point $\Lambda_\ast$.
        Without loss of generality, taking a subsequence further if necessary, we assume $\lim_{k\in\calKo\to\infty}\Lambda_k=\Lambda_\ast$.  We then have
        \begin{equation}
            \paren{H(x^\ast, \Lambda_\ast)}_{ij} \allowbreak =
            2\trace{\Acali(\xast) \Xastdagger  \Acalj(\xast) \Lambda_\ast}=
            \lim_{k\in \calK \rightarrow \infty}  \widehat{H}(\xk, \Lambda_k)_{ij}
            \label{eq:Hhatast}
        \end{equation}
        where
        $\hat{H}(\xk, \Lambda_k)$ and $H$ are defined in \eqref{eq:Hhat} and \cref{eq:sigma-term}, respectively.
\par
    To begin with, we assume $\Ker(X(\xast))\neq \{0\}$.
    From~\cref{eq:HessofLandpsimu}, we have
    \begin{align*}
         & \nablaij{L\paren{\xk, \allowbreak \Lambda_k}}
        \allowbreak                                        \\
         & = \paren{\nabla^2_{xx} \psimunu{\xk, \Zk}}_{ij}
        -  \allowbreak (1+\allowbreak\hnu(\muk)) \allowbreak \mu_k \allowbreak  \trace{\Acali\paren{\xk}
            \invXk\Acalj \paren{\xk}\invXk}.
    \end{align*}
    Choose $\ubar \in \calL(\xast)$ arbitrarily.
    From \cref{prop:innersemi}, there exists a sequence $\sbra{\uk}$ such that $\lim_{k\in\calKo\to\infty}\uk=\ubar$ and $\uk \in \Lastk$ for each $k\in \calKo$.
    Note that
    \begin{align}
        \begin{split}
            &\quadf{\HessLagrange{\xk, \Lambda_k}}{\uk}+\quadf{\widehat{H}(\xk, \Lambda_k)}{\uk} \\
            &=\quadf{\hesspsimu{\xk, \Zmu}}{\uk}
            +\underbrace{2 \trace{\Dx{\xk}{\uk}\invXkP \Dx{\xk}{\uk}\Lambda_k}}_{(A)}\\
            &\quad - \underbrace{\mukpn \trace{\Dx{\xk}{\uk}
                \invXkN
                \Dx{\xk}{\uk}
                \invXkN}}_{(B)} \\
            &\quad -\mukpn \trace{\Dx{\xk}{\uk}
                \invXkP
                \Dx{\xk}{\uk}
                \invXkP} \\
            &\quad -2\mukpn \trace{\Dx{\xk}{\uk}
                \invXkP\Dx{\xk}{\uk}
                \invXkN},
            \label{eq:traceexpand}
        \end{split}
    \end{align}
    where the equality follows from definition~\cref{eq:Lambda} of $\Lambda_k$ and \cref{eq:invDecomp}. Since $\uk \in \Lastk$, we have $\transpose{\Vkz}\Dx{\xk}{\uk}\Vkz = O$ which together with \eqref{def:invXkN} implies
    $(B) = 0$.
    Substituting definition\,\eqref{eq:Lambda} of $\Lambda_k$ into $(A)$ and using~\cref{eq:invDecomp} again imply
    \begin{equation*}
        \begin{aligned}\label{trace-3}
            (A)
             & = 2(1+\hnu(\muk))\muk\trace{\Dx{\xk}{\uk}
                \invXkP \Dx{\xk}{\uk}\invXkP}   \\
             & \quad +2\mukpn\trace{\Dx{\xk}{\uk}
                \invXkP \Dx{\xk}{\uk}\invXkN}   \\
             & \quad -2\hnu(\muk) \trace{\Dx{\xk}{\uk}
                \invXkP
                \Dx{\xk}{\uk}\Zk},
            \nonumber
        \end{aligned}
    \end{equation*}
    which together with~\cref{eq:traceexpand} and $(B)=0$ yields
    $$\quadf{\HessLagrange{\xk, \Lambda_k}}{\uk} + \quadf{\widehat{H}(\xk, \Lambda_k)}{\uk} = L_k + M_k + N_k,$$
    where, for each $k$,
    \begin{align*}
        L_k & \coloneqq  \quadf{\hesspsimu{\xk, \Zmu}}{\uk}, \nonumber                       \\
        M_k & \coloneqq \mukpn\trace{\Dx{\xk}{\uk} \invXkP \Dx{\xk}{\uk}\invXkP},\nonumber   \\
        N_k & \coloneqq -2\hnu(\muk) \trace{\Dx{\xk}{\uk} \invXkP \Dx{\xk}{\uk}\Zk}.\nonumber
    \end{align*}
    Hereafter, we aim to prove
    \begin{gather*}
        (\mathrm{F.}1)\; \liminf_{k \rightarrow \infty} L_k \geq 0,\; (\mathrm{F.}2)\; \liminf_{k\rightarrow \infty} M_k \geq 0,\; (\mathrm{F.}3)\;\lim_{k\rightarrow \infty} N_k = 0,
    \end{gather*}
    which together with \eqref{eq:Hhatast} actually entails $\quadf{\HessLagrange{\xast, \Lambda_\ast}}{\ubar} + \quadf{H(\xast, \Lambda_\ast)}{\ubar} \geq 0.$ Since $\bar{u}$ is an arbitrary point in $\mathcal{L}(x^\ast)$ and the KKT conditions hold at $\paren{\xast, \Lambda_\ast}$ from~\cref{theorem:KKT}, $\paren{\xast, \Lambda_\ast}$ is a {\wSOSP}. Thus, the proof is complete.
    Let us move on to the proofs of (F.1), (F.2), and (F.3).
    \begin{proof}[Proof of $(\mathrm{F.1})$]
        From condition~\eqref{eq:outsecond}, we have $
            L_k \geq -\heps_H(\mu_k) \allowbreak \bigl(1+\mu_k \fnorm{\invXk} + \fnorm{\Zk}\bigr)^2 \allowbreak \norm{\uk}^2,$ which together with~\eqref{eq:parameter}, $\lim_{k\to\infty}\mu_k=0$, and the boundedness of $\sbra{\uk}, \sbra{\Zk},$ and $\sbra{\muk \invXk}$ gives (F.1).
    \end{proof}
    \begin{proof}[Proof of $(\mathrm{F.2})$]
        Since $\invDkp =\paren{\paren{\Dkp}^{-1/2}}^2$, we obtain
        $$M_k\allowbreak=\allowbreak \mukpn\allowbreak \biggl\| \Bigl(\paren{\Dkp}^{-1/2} \allowbreak \transpose{\Vkp}\Bigr)\allowbreak \Dx{\xk}{\uk}\paren{\Vkp \paren{\Dkp}^{-1/2}}\biggr\|_{\mathrm{F}}^2 \geq 0
        $$
        and therefore (F.2) holds.
    \end{proof}
    \begin{proof}[Proof of $(\mathrm{F.3})$]
        In the case of $N_k$, we have
        \begin{align}
            \begin{split}
                \abs{N_k}
                \leq 2\hnu(\muk)  \fnorm{\Dx{\xk}{\uk} \invXkP}\fnorm{\Dx{\xk}{\uk}\Zk}.
            \end{split}
            \label{eq:Nbounded}
        \end{align}
        From~\Cref{assumption:bounded}, the continuity of $\Acali$ for all $i \in \sbra{1,\ldots, n}$, and the boundedness of $\sbra{\uk}$ and $\sbra{\Zk}$, we have the boundedness of $\sbra{\Dx{\xk}{\uk}}$ and $\sbra{\Dx{\xk}{\uk}\Zk}$.
        Recall $\lim_{k\in\calKo\to\infty}\Dkp=\Dastp\in \calD^{\rast}_{++}$ from equation\,\eqref{eq:limitmany}.
            Hence, $\sbra{\invDkp}$ is bounded, and thus $\sbra{\invXkP}$ is also bounded.
        These facts ensure that $\sbra{\Dx{\xk}{\uk}\invXkP}$ is bounded. Therefore by driving $k\in\calKo\to\infty$,
        $\hnu(\muk) \rightarrow 0$ holds and hence we obtain (F.3).
    \end{proof}
    \par
    We next consider the case where $\Ker(X(\xast))=\{0\}$, i.e., $X(\xast)\in S^m_{++}$. Then, we have
    $\Lastk=\calL(\xast)=\euclid^n$ for any $k\in\calKo$ and also have $\Lambda_{\ast}=O$ by the complementarity condition.
    Moreover, note that $H(x^{\ast},\Lambda_{\ast})=O$ holds.
    Then, by almost the same argument as above ignoring the terms including $X_k^N$ and $U_k^N$, we can reach the desired conclusion again.
\par
    The proof is complete.
    \qed
\end{proof}

\section{Worst-case iteration complexity}
\label{chapter:inner}
In this section, we analyze the worst-case iteration complexity of~\hyperlink{alg:inner}{Fixed-NC-IPM} for $\varepsilon$-SOSP$(\mu,\nu)$.
Throughout this section, we suppose that $\feasi\neq \emptyset$ and the initial point $(x^1,Z_1)\in \feasi\times \symmetric^m_{++}$ can be found.
Recall notation\,\eqref{eq:xell}.

The step size rules in Updates\,\ref{alg:dual},\ref{alg:primal}, and \ref{alg:negative_curvature} yield
\begin{align}
    \norm{\xell-\xellp} \leq \frac{\mineig{\Xell}}{2L_0},\ \fnorm{\Zell - \Zellp} \leq \frac{\mineig{\Zell}}{2}\label{al:prop2}
\end{align}
for each $\ell\ge 1$, which can be verified with direct calculation in each Update. For example, in Update~2, since $\alpha_{\ell}\le \frac{\mineig{\Xell}}{2L_0\norm{d_{\xell}}}$ by definition, we have 
  \[
  \norm{\xell-\xellp}
  = \alpha_{\ell}\norm{d_{x^{\ell}}}\le
  \frac{\mineig{\Xell}}{2L_0}.
  \]
Relevant to inequalities\,\eqref{al:prop2}, we introduce the following sets for each $\ell\ge 1$:
\begin{align*}
    \calNellx & \coloneqq\left\{x\in \euclid^n\middle| \norm{x-\xell} \leq \frac{\mineig{\Xell}}{2L_0}\right\},      \\
    \calNellZ & \coloneqq \left\{Z\in \symmetric^m \middle| \fnorm{Z-\Zell} \leq \frac{\mineig{\Zell}}{2}  \right\}.
\end{align*}
It is clear that for each $\ell\ge 1$
\begin{equation}
\xellp\in \calNellx,\ \Zellp\in\calNellZ \label{eq:calNellxZ}
\end{equation}
  from \eqref{al:prop2}. Hereafter, for the sake of brevity, we often write
\begin{equation}
X \coloneqq X(x).\label{eq:X(x)}
\end{equation}
\par
We next show that $\paren{\xell, \Zell}\in \feasi \times \symmetric^m_{++}$ holds for each $\ell\ge 1$.
\begin{proposition}
    \leavevmode
    \label{proposition:forinner}
    Suppose that \cref{assumption:Lipschtiz} holds. For each $\ell\ge 1$, the following properties hold.
    \begin{enumprop}
        \item 
        For any $x\in\calNellx$, it holds that
        $0<\frac{1}{2}\eigi{\Xell} \leq \eigi{\Xx} \leq \frac{3}{2} \eigi{\Xell}$ for each $i \in \sbra{1,\ldots,m}$.
        Therefore $\Xx \in \symmetric^{m}_{++}$, namely $x\in \feasi$.
        Besides, $\frac{1}{2}\fnorm{\Xell} \allowbreak\leq \fnorm{\Xx} \leq \frac{3}{2}\fnorm{\Xell}$ and $\fnorm{\invXx} \leq 2 \fnorm{\invXell}$ hold.\label{lemma:feasiblity}
        \item 
        For any $Z\in\calNellZ$, $\fnorm{\Zz^{-1}} \leq 2\fnorm{\Zell^{-1}}$
        and $\Zz \in \symmetric^m_{++}$ hold. \label{lemma:Zpd}
    \end{enumprop}
    In particular,
    $\emptyset \neq \calNellx\times \calNellZ\subseteq \feasi \times \symmetric^m_{++}$ for each $\ell$, and thus by \eqref{eq:calNellxZ} $\left\{\paren{\xell, \Zell}\right\} \allowbreak \subseteq  \feasi \times \symmetric^m_{++}$.
\end{proposition}
\begin{proof}
We first show \cref{lemma:feasiblity}. Note \eqref{eq:X(x)}.
From the Wielandt-Hoffman theorem~\cite[Corollary 6.3.8]{Horn2012} and~\cref{eq:Xlip}, we have, for each $i\in \{1,\ldots, m\}$, \allowbreak $\lvert\eigi{\Xell} \allowbreak -\allowbreak \eigi{\Xx}\rvert\allowbreak
        \leq \fnorm{\Xell-\Xx}\leq L_0 \norm{\xell - x}$,
    which together with the first inequality in \eqref{al:prop2}
    implies \[\abs{\eigi{\Xell} - \eigi{\Xx}}\allowbreak \leq \frac{\mineig{\Xell}}{2} \allowbreak
        \leq \frac{\eigi{\Xell}}{2},\]
    which leads to $\frac{\eigi{\Xell}}{2}\leq \eigi{\Xx}
        \leq  \frac{3\eigi{\Xell}}{2}.$
    We hence obtain $\frac{1}{2}\fnorm{\Xell} \leq \fnorm{\Xx} \leq \frac{3}{2}\fnorm{\Xell}$
    and $
        \fnorm{\invXx}
        \allowbreak =   \allowbreak \sqrt{\sum_{i=1}^n \paren{\frac{1}{\eigi{\Xx}}}^2}
        \allowbreak \leq  \allowbreak
        \sqrt{\sum_{i=1}^n \paren{\frac{2}{\eigi{\Xell}}}^2}
        \allowbreak  =   \allowbreak  2 \fnorm{\invXell}.$
        The proof of \Cref{lemma:Zpd} is omitted because it can be done in a similar manner to \cref{lemma:feasiblity}.
        \par
        Finally, we consider the last assertion.
    Since $(x^1,Z_1)\in \feasi\times \symmetric^m_{++}$, it follows that
    $\mathcal{N}_1(x^1)\times \mathcal{N}_2(Z_1)\neq \emptyset$.
    By \cref{lemma:feasiblity} and \cref{lemma:Zpd} with $\ell=1$, 
 $\emptyset \neq \calNellx\times\calNellZ \subseteq \feasi\times \symmetric^m_{++}$ is obtained for each $\ell$ inductively.
    Hence, we obtain the desired conclusion.
    \par
    The proof is complete.    \qed
\end{proof}
Later on, we will prove that if one of the approximate stationary conditions~\eqref{eq:aprox} is violated, there exists a direction that provides a sufficient decrease of the primal-dual merit function $\psi_{\mu,\nu}$. In order to prove this, we show the following ``local'' \Lipschitz continuities of the derivatives and the Hessian of $\psi_{\mu,\nu}$. We provide the proof of this lemma in the Appendix.
\begin{lemma}[local \Lipschitz continuities]
    \label{lemmabreak:localLipschitz}
    Suppose that \cref{assumption:Lipschtiz} holds. 
    For each $\ell\ge 1$, we have the following properties, where
    $l_Z(\tblue{\Zell})$, $l_x(\tblue{\xell,\Zell})$, and $l_{xx}(\tblue{\xell,\Zell})$ are defined in \autoref{alg:dual}, \autoref{alg:primal}, and \autoref{alg:negative_curvature}, respectively.
    \begin{enumlemma}
        \item \textbf{\textup{local \Lipschitz continuity of $X^{-1}$:}}
        \label{lemma:invX}
        For any $x\in\calNellx$,
        $$ \fnorm{\invXell-\invXx} \leq  \allowbreak 2 L_0 \allowbreak \fnorm{\invXell}^2 \allowbreak \norm{\xell - x}.
        $$
        \item \textbf{\textup{local \Lipschitz continuity of $\nabla_x \psi_{\mu,\nu}$:}} \label{lemma:nabalxLips}
        For any $x\in\calNellx$,
        \[
            \norm{\nabla_x \psimunu{\xell, \Zell} -\nabla_x \psimunu{x, \Zell}}
            \leq l_x(\xell,\Zell) \norm{\xell - x}.
        \]
        \item \textbf{\textup{local \Lipschitz continuity of $\nabla_Z \psi_{\mu,\nu}$:}} \label{lemma:nabalZLips}
        For any $Z\in\calNellZ$,
        \[
            \fnorm{\nabla_Z \psimunu{\xell, \Zell} - \nabla_Z \psimunu{\xell, \Zz}} \leq l_Z(\Zell)\fnorm{\Zell- \Zz}.
        \]
        \item \textbf{\textup{local \Lipschitz continuity of $\nabla^2_{xx} \psi_{\mu,\nu}$:}}  \label{lemma:hessLips}
        For any $x\in\calNellx$,
        \[\norm{\nabla^2_{xx} \psimunu{\xell, \Zell} - \nabla^2_{xx} \psimunu{x, \Zell}} \leq
            l_{xx}(\xell,\Zell)
            \norm{\xell - x}.\]
    \end{enumlemma}
\end{lemma}
\label{section:descent}
In order to derive a worst-case iteration complexity, we moreover show the following descent lemmas for Updates~\ref{alg:dual},~\ref{alg:primal}, and~\ref{alg:negative_curvature} by using~\cref{lemmabreak:localLipschitz}.
\begin{lemma}[descent lemmas]
    Suppose that \cref{assumption:Lipschtiz} holds. Then we have the following descent lemmas.
    \label{lemmabreak:descent}
    \begin{enumlemma}
        \item \textbf{\textup{descent lemma for~\autoref{alg:dual}:}} \label{descent:dual}
          Suppose that the approximated \comple condition ~\eqref{eq:apcomple} is violated, that is,
          \begin{equation}
            \fnorm{\nabla_Z \psimunu{\xell,\Zell}}\allowbreak >
            \varepsilon_\mu (1+\mu \fnorm{{\Zell}^{-1}}). \label{eq:nablaZpsi}
\end{equation}
          Then,~\autoref{alg:dual} decreases $\psi_{\mu,\nu}$ at least by
          $$\sigma_1 \coloneqq
            \min \Bigl\{\frac{\mu \varepsilon_\mu\iotamin}{4\iotamax},\allowbreak \frac{\mu \epsmu^2 \iotamin^2}{4\nu\iotamax^2}\Bigr\}.$$
        \item \textbf{\textup{descent lemma for~\autoref{alg:primal}:}} \label{descent:primal}
        Suppose that the approximated \statio condition ~\eqref{eq:apsta} is violated, that is,
\begin{equation}
  \norm{\nabla_x \psimunu{\xell, \Zell}}\allowbreak>  \varepsilon_g\bigl(1+\allowbreak\mu\fnorm{\invXell}\allowbreak+\allowbreak\fnorm{\Zell}\bigr).\label{eq:nablaxpsimu}
\end{equation}
  Then,~\autoref{alg:primal} decreases $\psi_{\mu,\nu}$ at least by
        \begin{align*}
            \sigma_2\coloneqq
            \min
            \biggl\{\frac{\mu \varepsilon_g \hmin}{4 L_0 \hmax}, & \frac{\varepsilon_g^2 \hmin^2}{8L_1 \hmax^2}, \frac{\varepsilon_g^2 \hmin^2}{4\nu L_1 \hmax^2}, \\
                                                                 & \frac{\mu \varepsilon_g^2 \hmin^2}{16(1+\nu)L_0^2 \hmax^2},
            \frac{\varepsilon_g^2 \hmin^2}{4\paren{1+\nu}L_1 \hmax^2 }\biggr\}.
        \end{align*}
        \item \textbf{\textup{descent lemma for~\autoref{alg:negative_curvature}:}} \label{descent:negative}
          Suppose that the approximated \secondorder condition~\eqref{eq:apsecond} is violated, that is,
          \begin{equation}
            \mineig{\hesspsimu{\xell, \Zell}}\allowbreak < \allowbreak-\varepsilon_H \allowbreak \bigl(1+\allowbreak \mu \fnorm{\invXell} \allowbreak+\allowbreak \fnorm{\Zell}\bigr)^2.\label{eq:hesspsi}
\end{equation}
            Then,~\autoref{alg:negative_curvature} decreases $\psi_{\mu,\nu}$ at least by
        \begin{align*}
            \sigma_3 \coloneqq
            \min \Biggl\{ & \frac{\mu^2 \varepsilon_H}{24L_0^2},
            \frac{2\varepsilon_H^3}{75L_2^2},
            \\
                          & \quad \frac{2\varepsilon_H^3}{5\nu^2 L_2^2},
            \frac{2\varepsilon_H^3}{5(1+\nu)^2 L_2^2},
            \frac{\mu^2\varepsilon_H^3}{40(1+\nu)^2 L_1^2 L_0^2},
            \frac{\mu^4\varepsilon_H^3}{1350(1+\nu)^2 L_0^6}\Biggr\}.
        \end{align*}
    \end{enumlemma}
\end{lemma}
\begin{proof}[Proof of \cref{descent:dual} ]
    The first-order descent lemma~\cite[Lemma 5.7]{Beck2017} together with \cref{lemma:nabalZLips} yields, for any $Z\in\calNellZ$,
        \[
            \psimunu{\xell, Z} \leq \psimunu{\xell, \Zell}
            + \innerstepsize \nabla_Z \psimunu{\xell, \Zell}^\top \innerZdirection
            + \frac{\innerstepsize^2 l_Z(\Zell)}{2}\norm{\innerZdirection}^2,
        \] which, together with
    $\innerZdirection = -\frac{\iotamin}{\iotamax^2}\calHell{\nabla_Z \psimunu{\Zell}}$ and $Z=\Zellp$, implies
    \begin{align*}
        & \psimunu{\xell, \Zellp} \\
        &\leq \psimunu{\xell, \tblue{\Zell}}           
                               - \frac{\innerstepsize \iotamin}{\iotamax^2}\left \langle \nabla_Z \psimunu{\xell, \Zell},  \calHell{\nabla_Z \psimunu{\xell, \Zell}}\right \rangle \\
                                &\quad 
        + \frac{\innerstepsize^2 \iotamin^2 l_Z(\Zell)}{2 \iotamax^4}\fnorm{\calHell{\nabla_Z \psimunu{\xell, \Zell}}}^2                                                              \\
                                & \leq  \psimunu{\xell, \tblue{\Zell}}                         -\frac{\innerstepsize \iotamin^2}{\iotamax^2}\fnorm{\nabla_Z \psimunu{\tblue{\xell}, \eta_{\ell}}}
        \paren{1- \frac{l_Z(\Zell) \innerstepsize}{2}},
    \end{align*}
    where the last inequality follows from \cref{eq:calHell} and the symmetry of $\calHell{\cdot}$. Since $\innerstepsize \leq \frac{1}{l_Z(\Zell)}$ by the step size rule,
    we have $
    \psimunu{\xell, \Zellp} \leq \psimunu{\xell, \Zell} \allowbreak- \frac{\innerstepsize \iotamin^2}{2 \iotamax^2} \fnorm{\nabla_Z \psimunu{\xell, \Zell}}^2.$ Therefore we see that  $\psi_{\mu,\nu}$ decreases at least by
    $$
    \xi_{\ell} \coloneqq \frac{\alpha_{\ell} \iotamin^2}{2 \iotamax^2}\fnorm{\nabla_Z \psimunu{\xl, \eta_{\ell}}}.$$
    In what follows, we evaluate a lower bound of $\xi_{\ell}$ for each choice of $\innerstepsize=\frac{\mineig{\Zell}}{2 \fnorm{\innerZdirection}},\ \frac{1}{l_Z(\Zell)}.$
    
    We first consider the choice of $\innerstepsize = \frac{\mineig{\Zell}}{2 \fnorm{\innerZdirection}}$, which together with
    $
        \fnorm{\innerZdirection } \allowbreak= \allowbreak \frac{\iotamin}{\iotamax^2}\fnorm{\calHell{\nabla_Z \psimunu{\xell, \Zell}}}
        \allowbreak \leq \allowbreak \frac{\iotamin}{\iotamax} \fnorm{\nabla_Z \psimunu{\xell, \Zell}}
    $ and \eqref{eq:calHell} implies
    \begin{align}
        \begin{split}
            \xi_{\ell}
            &=\frac{\mineig{\Zell} \iotamin^2}{4 \fnorm{\innerZdirection} \iotamax^2}
            \fnorm{\nabla_Z \psimunu{\xell, \Zell}}^2\\
            &> \frac{\mineig{\Zell} \iotamin}{4  \iotamax} \varepsilon_\mu (1+\mu \fnorm{\Zell^{-1}}) > \frac{\mu \varepsilon_\mu \mineig{\Zell}\fnorm{\Zell^{-1}}\iotamin }{4 \iotamax}
            \geq \frac{\mu \varepsilon_\mu \iotamin}{4 \iotamax},
            \nonumber
        \end{split}
    \end{align}
    where the first inequality follows from \eqref{eq:nablaZpsi}.

 Next, we consider the choice of  $\alpha_{\ell} = \frac{1}{l_Z(\Zell)}= \frac{1}{2\mu \nu \fnorm{\Zell^{-1}}^2}$, which combined with \eqref{eq:nablaZpsi} yields
    \begin{align*}
        \begin{split}
            \xi_{\ell}
            &=\frac{\iotamin^2}{4 \mu\nu \fnorm{\Zell^{-1}}^2 \iotamax^2}\fnorm{\nabla_Z \psimunu{\xell, \Zell}}^2 \\
            &>\frac{\iotamin^2}{4 \mu\nu \fnorm{\Zell^{-1}}^2 \iotamax^2} \epsmu^2 (1+\mu \fnorm{\Zell^{-1}})^2 >  \frac{\mu \epsmu^2 \iotamin^2}{4\nu \iotamax^2}.
        \end{split}
    \end{align*}
    These lead to the desired conclusion and the proof is complete.
\end{proof}
\begin{proof}[Proof of \cref{descent:primal}]
    The first-order descent lemma~\cite[Lemma 5.7]{Beck2017} with \cref{lemma:nabalxLips} yields that  
        \[
            \psi_{\mu,\nu}\bigl(x,\allowbreak \Zell\bigr)\allowbreak \leq \allowbreak \psi_{\mu,\nu} \bigl(\xellp, \allowbreak \Zell\bigr)
            \allowbreak+\allowbreak \innerstepsize \nabla_x \psimunu{\xell, \Zell}^\top \innerdirection
            + \frac{\innerstepsize^2 l_x(\xell, \Zell)}{2}\allowbreak\norm{\innerdirection}^2
            \]
            for any $x\in\calNellx$.
            Since $\xellp\in \calNellx$ by \eqref{eq:calNellxZ}, $\xellp$ can be set as $x$ above. Furthermore, substituting
              \begin{equation}
                \innerdirection = -\frac{\hmin}{\hmax^2}\Hell \allowbreak \nabla_x \psi_{\mu,\nu} \bigl(\xellp \allowbreak, \allowbreak \Zell\bigr), \label{eq:innerdir}
                \end{equation}
                we obtain
        \begin{align*}
            &\psimunu{\xellp, \Zell} \\
            & \leq \psimunu{\xellp, \Zell}
            - \frac{\innerstepsize \hmin}{\hmax^2}\nabla_x \psimunu{\xell, \Zell}^\top\Hell\nabla_x \psimunu{\xell, \Zell}                                                       \\
                                    &\quad  +\frac{\innerstepsize^2 \hmin^2 l_x(\xell, \Zell)}{2\hmax^4} \nabla_x \psimunu{\xell, \Zell}^\top \Hell^{2}\nabla_x \psimunu{\xell, \Zell} \\
                                    & \leq  \psimunu{\xellp, \Zell}                      -
            \frac{\innerstepsize \hmin^2}{\hmax^2}
            \norm{\nabla_x \psimunu{\xell, \Zell}}^2
            \paren{1-\frac{l_x(\xell, \Zell) \innerstepsize }{2}},
        \end{align*}
        where the last inequality follows from \cref{eq:Hell}.
    Since $\innerstepsize \leq \frac{1}{l_x(\xell, \Zell)}$ follows from the step size rule,
    the above inequality further implies
    $$
    \psi_{\mu,\nu}\bigl(\xellp,\allowbreak \Zell\bigr) \allowbreak \leq \allowbreak \psimunu{\xell, \Zell}\allowbreak -\allowbreak \frac{\innerstepsize \hmin^2}{2 \hmax^2}\allowbreak \norm{\nabla_x \psimunu{\xell, \Zell}}^2.$$
    Now, we evaluate a lower bound of the function decrease
    $$\delta_{\ell}\coloneqq
    \allowbreak \frac{\innerstepsize \hmin^2}{2 \hmax^2} \norm{\nabla_x \psimunu{\xell, \Zell}}^2$$
    for each choice of $\innerstepsize= \frac{\mineig{\Xell}}{2 L_0\nnorm{\innerdirection}}$ and $ \frac{1}{l_x \paren{\xell, \Zell}}.$
    We first consider the choice of $\innerstepsize = \frac{\mineig{\Xell}}
    {2 L_0 \nnorm{\innerdirection}}$, which together with
    \eqref{eq:nablaxpsimu} and
    $
\norm{\innerdirection}\leq \frac{\hmin}{\hmax}\norm{\nabla_x \psimunu{\xell, \Zell}}
    $ from \eqref{eq:innerdir} and \eqref{eq:Hell}
    implies
    \begin{align*}
            & \delta_{\ell} = \frac{\hmin^2}{2\hmax^2} \paren{\frac{\mineig{\Xell}}{2 L_0  \norm{\innerdirection}}}
            \norm{\nabla_x \psimunu{\xell, \Zell}}^2\\
            &\geq\frac{\varepsilon_g \hmin}{2 \hmax} \paren{\frac{\mineig{\Xell}}{2 L_0 }} \paren{1+\mu\fnorm{\invXell}+\fnorm{\Zell}}\\
            &\geq\frac{\varepsilon_g \mu \hmin }{4 L_0 \hmax}\mineig{\Xell}\fnorm{\invXell}\\
            &\geq\frac{\varepsilon_g \mu \hmin}{4 L_0 \hmax},
    \end{align*}
    where the last inequality is due to $\mineig{\Xell}^{-1}\le \fnorm{\invXell}.$

    Next, we consider the choice of $\innerstepsize = \frac{1}{l_x(\xell, \Zell)}$, which yields
    \begin{align*}
            \delta_{\ell}
            &> \frac{\hmin^2}{2 \hmax^2 l_x(\xell, \Zell)}
            \varepsilon^2_g\paren{1+\mu\fnorm{\invXell}+\fnorm{\Zell}}^2\\
            &=\frac{\hmin^2\varepsilon^2_g\paren{1+\mu\fnorm{\invXell}+\fnorm{\Zell}}^2}
            {2 \hmax^2\paren{L_1 + \nu L_1\fnorm{\Zell}+ 2\mupn L_0^2\fnorm{\invXell}^2+ \mupn L_1 \fnorm{\invXell}}}\\
            &\geq \frac{\hmin^2}{2 \hmax^2} \varepsilon^2_g\paren{1+\mu\fnorm{\invXell}+\fnorm{\Zell}}^2\\
            & \cdot \min
            \sbra{\frac{1}{4L_1}, \frac{1}{4\nu L_1 \fnorm{\Zell}}, \frac{1}{8\mupn L_0^2\fnorm{\invXell}^2},
                \frac{1}{4\mupn L_1\fnorm{\invXell} }}\\
            &\geq\min
            \sbra{\frac{\varepsilon_g^2 \hmin^2}{8 L_1 \hmax^2},
                \frac{\varepsilon_g^2 \hmin^2}{4\nu L_1 \hmax^2},
                \frac{\mu \varepsilon_g^2 \hmin^2}{16(1+\nu) L_0^2 \hmax^2},
                \frac{\varepsilon_g^2 \hmin^2}{4\paren{1+\nu} L_1 \hmax^2 }}
            \nonumber,
    \end{align*}
    where the second inequality is derived from $\frac{1}{a+b+c+d} \geq \frac{1}{4} \min \sbra{\frac{1}{a}, \frac{1}{b}, \frac{1}{c}, \frac{1}{d}}$ for $a, b, c, d>0$,
    and the third one is from 
    \[
    \paren{1+\mu\fnorm{\invXell}+\fnorm{\Zell}}^2\ge \max\left(
    1,2\fnorm{\Zell},\mu^2\fnorm{\Xell^{-1}}^2,2\mu\fnorm{\Xell^{-1}}\right).
    \]
    These lead to the conclusion and the proof is complete.
\end{proof}
\begin{proof}[Proof of \cref{descent:negative}]
    Since $\innerdirection$ is a normalized eigenvector corresponding to
    $\lambda_{\rm min}\bigl(\nabla^2_{xx} \psi_{\mu,\nu}(\xell,\allowbreak \Zell)\bigr)$ such that $\innerdirection^\top \nabla_x \psi_{\mu,\nu} (\xell,\allowbreak \Zell) \leq 0$, the second-order descent lemma~\cite[Lemma 1]{Nesterov2006} together with \cref{lemma:hessLips} yields, for any $x\in\calNellx$,
    \begin{align}
        \psimunu{x, \Zell}
        \leq
        \psimunu{\xell, \Zell}+
        {\alpha^2_\ell}
        \paren{
            \frac{\mineig{\hesspsimu{\xell, \Zell}}}{2}
            +
            \frac{l_{xx}(\xell, \Zell) \innerstepsize}{6}}.
        \label{eq:eigend}
    \end{align}
    Note that we have 
    \begin{align}
         & \frac{\mineig{\hesspsimu{\xell, \Zell}}}{2}+
        \frac{l_{xx}(\xell, \Zell) \innerstepsize}{6}<0
        \label{eq:eigendecrease}
    \end{align}
    from \eqref{eq:hesspsi}, $l_{xx} (\xell, \Zell)>0$, and $0 < \innerstepsize \leq -\frac{2\mineig{\hesspsimu{\xell, \Zell}}}{l_{xx}(\xell, \Zell)}$.
      By using \cref{eq:hesspsi,eq:eigend,eq:eigendecrease}, $\psi_{\mu,\nu}$ decreases at least by
    \[
    \eta_{\ell} \coloneqq \allowbreak -\alpha^2_{\ell} \frac{\mineig{\hesspsimu{\xell, \Zell}}}{6} > 0.
    \]
    We evaluate a lower bound of $\eta_{\ell}$ for each choice of $\innerstepsize=-\frac{2 \mineig{\hesspsimu{\xell, \Zell}}}{l_{xx}(\xell, \Zell)}$ and $\frac{\mineig{\Xell}}{2L_0 \norm{\innerdirection}}.$
    \par
    We first consider the choice of  $\innerstepsize=   -\frac{2 \mineig{\hesspsimu{\xell, \Zell}}}{l_{xx}(\xell, \Zell)}$.
Recalling \eqref{eq:hesspsi} and the definition of $l_{xx}(\xell, \Zell)$ in Update~3, we have
    \begin{align*}
         & \eta_{\ell}
        =-\frac{2\paren{\mineig{\hesspsimu{\xell, \Zell}}}^3}{3 l_{xx}(\xell, \Zell)^2}     \\
         & \geq \frac{2\varepsilon_H^3 \paren{1+\mu \fnorm{\invXell} + \fnorm{\Zell}}^6}{3}
        \Bigl( L_2+
        \nu L_2 \fnorm{\Zell}
        \\
         & \quad +
        \mu \paren{1+\nu}\paren{ L_2 \fnorm{\invXell}+
            4 L_1 L_0\fnorm{\invXell}^2 + 6 L_0^3 \fnorm{\invXell}^3}\Bigr)^{-2}            \\
         & \geq \frac{2\varepsilon_H^3}{75}\paren{1+\mu \fnorm{\invXell} + \fnorm{\Zell}}^6 \\
         & \quad \cdot
        \min \left\{\frac{1}{L_2^2}, \frac{1}{\nu^2 L_2^2 \fnorm{\Zell}^2}, \frac{1}{\mu^2 (1+\nu)^2 L_2^2
            \fnorm{\invXell}^2},\right.                                                     \\
         & \quad \quad \quad \quad \left. \frac{1}{16\mu^2 (1+\nu)^2 L_1^2 L_0^2 \fnorm{\invXell}^4},
        \frac{1}{36 \mu^2 (1+\nu)^2 L_0^6 \fnorm{\invXell}^6}\right\}                       \\
         & \geq \frac{2\varepsilon_H^3}{75}
        \min \sbra{\frac{1}{L_2^2}, \frac{15}{\nu^2 L_2^2 },
            \frac{15}{(1+\nu)^2 L_2^2},
            \frac{15\mu^2}{16 (1+\nu)^2 L_1^2 L_0^2},
            \frac{\mu^4}{36  (1+\nu)^2 L_0^6}},
    \end{align*}
    where the second inequality follows from 
    $\frac{1}{a+b+c+d+e} \geq  \frac{1}{5} \min \sbra{\frac{1}{a}, \frac{1}{b}, \frac{1}{c}, \frac{1}{d}, \frac{1}{e}}$
      for $a,b, c, d, e > 0$, and the third one does from
      \[
    \paren{1+\mu \fnorm{\invXell} + \fnorm{\Zell}}^6\geq \max\left(
1, 15\fnorm{\Zell}^2,15\mu^2\fnorm{\Xell^{-1}}^2,15\mu^4\fnorm{\Xell^{-1}}^4,\mu^6\fnorm{\Xell^{-1}}^6
      \right).
      \]
      \par
    Next, we consider the choice of  $\innerstepsize =  \frac{\lambda_{\min}(\Xell)}{2L_0 \norm{\innerdirection}}$, which 
    together with $\norm{\innerdirection} = 1$ and \eqref{eq:hesspsi} implies
    \begin{align*}
        \eta_{\ell}
         & = -\paren{ \frac{\mineig{\Xell}}{2L_0}}^2 \frac{\mineig{\hesspsimu{\xell, \Zell}}}{6}                                 \\
         & \geq \paren{ \frac{\mineig{\Xell})}{2L_0}}^2 \frac{\varepsilon_H \paren{1+\mu \fnorm{\invXell} + \fnorm{\Zell}}^2}{6} \\
         & \geq \paren{ \frac{\mineig{\Xell}}{2L_0}}^2
        \frac{\varepsilon_H \paren{\mu^2 \fnorm{\invXell}^2}}{6}\\
        &\geq \frac{\mu^2 \varepsilon_H}{24L_0^2},
    \end{align*}
    where the last inequality follows from
    $\fnorm{\invXell}\ge \lambda_{\rm min}(\Xell)^{-1}$.
    These lead to the conclusion and the proof is complete. \qed
\end{proof}
To establish the convergence theorem for~\hyperlink{alg:inner}{Fixed-NC-IPM}, let us make one more assumption.
\begin{assumption}
    For each $\mu,\nu > 0$,  $\psi_{\mu,\nu}$ is bounded below on $\feasi \times \symmetric^m_{++}.$
    That is, with some $\psi_{\mu,\nu}^\ast$, we have $\psimunu{x, Z} \geq \psi_{\mu,\nu}^\ast$ for all $(x,Z) \in \feasi \times \symmetric^m_{++}$.
    \label{assumption:merit_bound}
\end{assumption}
\begin{theorem}[iteration complexity for approximate SOSPs]\label{thm:comp}
Suppose that~\cref{assumption:Lipschtiz,assumption:merit_bound} hold. Then,
  given positive parameters $\mu,\nu$, \hyperlink{alg:inner}{Fixed-NC-IPM} terminates at $\varepsilon$-SOSP$(\mu,\allowbreak \nu)$
    within $\ovzeta:=\frac{\psimunu{x_1, Z_1} - \psi_{\mu,\nu}^\ast}{\min \sbra{\sigma_1, \sigma_2, \sigma_3}}$ iterations, where
    $\varepsilon=(\varepsilon_g, \epsmu, \varepsilon_H)$, and
$\sigma_1,\sigma_2$ and $\sigma_3$ are defined in~\cref{lemmabreak:descent}. In particular, 
    we have $\ovzeta=\order{\frac{\psimunu{x_1, Z_1} - \psi_{\mu,\nu}^\ast}{\min\sbra{ \varepsilon_g^2,  \epsmu^2, \varepsilon_H^3}}}$, namely,
there exists a constant $M>0$ such that
$\ovzeta\le M
\left|\frac{\psimunu{x_1, Z_1} - \psi_{\mu,\nu}^\ast}{\min\sbra{ \varepsilon_g^2,  \epsmu^2, \varepsilon_H^3}}\right|$
for any $\varepsilon_g,\epsmu$, and $\varepsilon_H>0$ sufficiently small.
\end{theorem}
\begin{proof}
We first show the first assertion. By~\cref{lemmabreak:descent}, $\psimunu{\cdot, \cdot}$ decreases by at least $\min \sbra{\sigma_1, \sigma_2, \sigma_3}>0$ every iteration, as long as \cref{eq:apsta,eq:apcomple,eq:apsecond} are violated. On the other hand, due to \cref{assumption:merit_bound}, the total amount of the decrease is bounded above by $\psimunu{x_1, Z_1} - \psi_{\mu,\nu}^\ast$.     
If either one of \cref{eq:apsta,eq:apcomple,eq:apsecond} is violated after $\ovzeta$ iterations,
$\psimunu{\cdot, \cdot}$ must decreases further by at least $\min \sbra{\sigma_1, \sigma_2, \sigma_3}$, resulting in a value of $\psi_{\mu,\nu}$ smaller than $\psi_{\mu,\nu}^\ast$.
This \tteal{contradicts} \cref{assumption:merit_bound}. We thus establish the desired assertion.
\par
Next, we show the second assertion. Notice that 
\[
  \sigma_1 = \order{\epsmu^2},\ \sigma_2 = \order{\varepsilon_{g}^2},\ \sigma_3 = \order{\varepsilon_{H}^3}
\]
  follow since $\mu,\nu$ together with $h_{\{{\rm min},{\rm max}\}}, \kappa_{\{{\rm min},{\rm max}\}}$ are algorithmic parameters and $L_0,L_1,L_2$ are constants dependent on the problem data. Therefore, we obtain
    \[\frac{\psimunu{x_1, Z_1} - \psi_{\mu,\nu}^\ast}{\min \sbra{\sigma_1, \sigma_2, \sigma_3}} =  \order{\frac{\psimunu{x_1, Z_1} - \psi_{\mu,\nu}^\ast}{\min\sbra{ \varepsilon_g^2,  \epsmu^2, \varepsilon_H^3}}}.\]
    \par
    The proof is complete. \qed
\end{proof}

\begin{remark}
    When the objective $f$ is convex and the constraint $G$ is affine, Update~3 is not necessary since KKT points are SOSPs.
    Hence, in this case, $\sigma_3$ and $\varepsilon_H$ from Update~3 can be neglected, and then $\ovzeta$ reduces to
  $\frac{\psimunu{x_1, Z_1} - \psi_{\mu,\nu}^\ast}{\min \sbra{\sigma_1, \sigma_2}}$ or $\order{\frac{\psimunu{x_1, Z_1} - \psi_{\mu,\nu}^\ast}{\min\sbra{ \varepsilon_g^2,  \epsmu^2}}}.$
\end{remark}
\begin{remark}
    In practice, it is known that more iterations are required
for minimizing $\psi_{\mu,\nu}$ as $\mu$ gets smaller. Interestingly, this phenomenon can be induced from \cref{thm:comp}.
    Indeed, consider the number of iterations for finding
    $\varepsilon$-SOSP$(\mu,\allowbreak \nu)$, $\ovzeta$ from \cref{thm:comp}, while decreasing $\mu$ but fixing $\nu$ and $\varepsilon$. Since 
    $
  \sigma_1 = \order{\mu},\ \sigma_2 = \order{\mu}$,\ and $\sigma_3 = \order{\mu^4}
  $ hold 
  by the definition of $\sigma_1,\sigma_2,\sigma_3$, we can lower-bound
  $\ovzeta$ with $\order{\mu^{-4}(\psimunu{x_1, Z_1} - \psi_{\mu,\nu}^\ast)}$, which implies\footnote{
    We suppose that 
$\psimunu{x_1, Z_1} - \psi_{\mu,\nu}^\ast$ does not converge to 0 as $\mu\to 0$.
} that the number of iterations for approximate SOSPs increases as $\mu$ approaches 0. 
\end{remark}

\section{Extension to primal IPM}
\label{sec:primal}
We can derive a primal IPM from
the proposed primal-dual IPM with the following modification.
In \cref{alg:outer}, we fix $\hnu(\muk)\equiv 0$ for every $k$, and thus $\psi_{\mu,0}$ is the merit function under this setting.
In \cref{alg:inner}, we set $\Zell \coloneqq \mu \invXell$ at each iteration so that we have $\nabla_Z\psi_{\mu,0}(\xell,\Zell)=O$ for every $\ell$ (see \eqref{eq:gradfZ})  and remove Procedure~1.
\Cref{alg:innerprimal} describes the resulting algorithm.
\par
With these changes, it is possible to establish
  the convergence of the primal IPM to SOSPs with an iteration complexity in almost the same argument as the primal-dual IPM.
Specifically, 
Fixed-NC-PIPM (\cref{alg:innerprimal}) converges to an $\varepsilon$-SOSP$(\mu,0)$ within $\order{\frac{\psi_{\mu,0}\paren{x_1, Z_1} - \psi_{\mu,0}^\ast}{\min\sbra{ \varepsilon_g^2,  \varepsilon_H^3}}}$ iterations under \cref{assumption:Lipschtiz} and \cref{assumption:merit_bound} with $\nu=0$. Since $\nabla_Z\psi_{\mu,0}(\xell,\Zell)=O$ for every $\ell$, the term $\epsmu^2$ is no longer necessary in \cref{thm:comp}.
  As well, the convergence of $\varepsilon$-SOSP$(\mu, 0)$ to {\wSOSP} can be proved under the same assumptions as
  \cref{theorem:SOSP}, but without \cref{assumption:zmuinvxbounded} and \cref{assumption:compleconv}.
  Indeed, $Z_k=\mu_k\Xk^{-1}$ unconditionally yields \cref{assumption:zmuinvxbounded} together with the assertion in \cref{theorem:Zlambda}.  
As a result, \cref{assumption:compleconv} becomes redundant
because it is required by the primal-dual IPM only for the sake of establishing \cref{theorem:Zlambda}.
\begin{algorithm}[htbp]
  \DontPrintSemicolon
  \caption{Fixed Negative Curvature \textbf{Primal} Interior Point Method with $\mu$, $\nu$, $\varepsilon_g$,
    $\varepsilon_\mu$, $\varepsilon_H$ (Fixed-NC-PIPM ($\mu,\nu, \varepsilon_g, \varepsilon_\mu,\varepsilon_H$))}
  \label{alg:innerprimal}
  \KwIn{$\xinit, Z_{1} = \mu X(\xinit)^{-1}, \mu, \nu, \varepsilon_g, \varepsilon_\mu, \varepsilon_H, L_0, L_1, L_2, \hmin, \hmax, \iotamin, \iotamax$}\tcp{$X(\xinit), Z_1\in \mathbb{S}^m_{++}$}
  \For{$\ell = 1, \ldots$ }{
    \uIf{$\norm{\nabla_x \psimunutwo{\xell, \Zell}}>   \varepsilon_g\paren{1+\mu\fnorm{\invXell}+\fnorm{\Zell}}$\tcp{Procedure 2}}{
      Set $\xellp$ by Update 2\;
      $\Zellp \leftarrow \mu \invXellp$
    }
    \uElseIf{$\mineig{\hesspsimutwo{\xell, \Zell}} < -\varepsilon_H \paren{1+\mu \fnorm{\invXell} + \fnorm{\Zell}}^2$\;\tcp{Procedure 3}}{
      Set $\xellp$ by Update~3\;
      $\Zellp \leftarrow \mu \invXellp$
    }
    \Else{
      Output$(\xell, \Zell)$
    }
  }
\end{algorithm}

\section{Numerical Experiments}
\label{section:numexp}
We conduct numerical experiments to examine the efficiency of directions of negative curvature.
For the sake of practical implementation, we integrate
line-search procedures: We use the upper bound of the step sizes $\innerstepsize$ given in \cref{proposition:forinner} as initial values and backtrack until $\pmunu$ decreases by 
\begin{itemize}[label=\textbullet]
\item  $\frac{\innerstepsize \hmin^2}{2 \hmax^2} \allowbreak \norm{\nabla_x \psimunu{\xell,\Zell}}^2$ for the scaled
gradient w.r.t. $x$,
\item $\frac{\innerstepsize \iotamin^2}{2 \iotamax^2} \allowbreak \bigl\| \nabla_Z \psi_{\mu,\nu} (\xell, \allowbreak \Zell)\bigr\|_{\mathrm{F}}^2$
for the scaled  gradient w.r.t. $Z$,
\item  $\innerstepsize^2 \frac{\mineig{\hesspsimu{\xell, \Zell}}}{6}$ for the negative curvature.
\end{itemize}
All experiments are implemented in Python with JAX~\cite{jax2018github} on a computer with AMD Ryzen 7 2700X and Ubuntu 18.04.
\subsection{Problem setting: shifted positive semidefinite factorization}
We consider positive semidefinite factorization (PSF) problems~\cite{Vandaele2018}, which are \tteal{a} generalization of nonnegative matrix factorization (NMF) problems~\cite{Paatero1994,Lee1999}. As well as NMF problems, PSF problems are nonconvex. Let $V$ be an $m\times n$ nonnegative matrix and $q < \min \sbra{m,n}$.
We solve the following PSF problem shifted by adding $r I_q$ to the constraints:
\begin{mini}
    {A_i, B_j}{\sum_{i=1}^m \sum_{j=1}^n \paren{V_{ij}-\langle A_i, B_j \rangle}^2}
    {\label{problem:psf_shfted}}{}
    \addConstraint{A_i + r I_q}{\in \symmetric^q_{+},\quad}{i = 1,\ldots,m}
    \addConstraint{B_j + r I_q}{\in \symmetric^q_{+},\quad }{j = 1,\ldots, n},
\end{mini}
where $r$ is a positive constant and $\{A_i\},\{B_j\}$ are matrices to be determined.
The origin, an interior point,  
is a strict saddle point of problem~\cref{problem:psf_shfted} unless $V=O$.
\par
To produce the matrix $V \in \euclid^{m\times n}$, we repeat the following procedures in order until the nonnegativity of $V$ is ensured.
\begin{enumerate}[label={(\arabic*)}]
    \item Generate $A^\ast_i \in \symmetric^q$ and $B^\ast_j \in \symmetric^q$ by multiplying a matrix with elements randomly generated from the uniform distribution over $[0, 1)$ and its transpose.
    \item Randomly set $20\%$ of eigenvalues of $A^\ast_i$ and $B^\ast_j$ as $-r$ so that the solution of problem~\eqref{problem:psf_shfted} is on the boundary.
    \item  Generate $V$ by
          $ V_{ij} \coloneqq \langle A^\ast_i, B^\ast_j \rangle$ for each $i,j \in \sbra{1, \ldots ,m} \times \sbra{1,\ldots , n}$.	
\end{enumerate}
\subsection{Details of the algorithms}
\label{subsection:numexp}
In order to examine the effects of directions of negative curvature, we also implement \hyperref[alg:outer]{Decreasing-NC-PDIPM} without negative curvature, namely, Procedure 3 in \hyperref[alg:inner]{Fixed-NC-PDIPM}.
We set $L_0 = 1$ and also set an initial value of $\mu$ to $0.3$. With respect to scaling of gradient directions,  we set $\Hell = I_n$ and $\calHell{\cdot}$ to be an identity mapping.
The update rules for the parameters are as follows:
\begin{gather*}
    \mukp := \min (0.8 \muk, 10  \muk ^{1.5}),
    \hnu (\mu) := \mu^{0.1},
    \heps_g(\mu) := \mu, \\
    \heps_H (\mu):=\mu,
    \heps_\mu (\mu) := \mu^{1.2}.
\end{gather*}
\par
Let $x\in \euclid^{\binom{q+1}{2}(n+m)}$ be the vector obtained by combining the vectorization of $A_i \allowbreak (i=1,\allowbreak \ldots,n)$ and $B_j (j=1,\ldots,m)$.
We have different ways of choosing the initial point $x^1$ in the two experimental settings.
In~\Cref{subsubsection:saddle_neighbor}
, we choose an initial point $x^1\in \euclid^{\binom{q+1}{2}(n+m)}$ randomly from a uniform distribution over $[0, 10^{-6})$ in order to place $x^1$ in the neighborhood of the strict saddle point so that there will be the use of directions of negative curvature.
From the above generation of $x^1$, $x^1$ is not guaranteed to be feasible. Therefore we generate $x^1$ until we find feasible $x^1$.
In~\Cref{subsubsection:comparison_with_pgd}, we choose an initial point $x^1 = O$ so that algorithms start from the strict saddle point.
\par 
With regard to $Z$ in~\cref{alg:outer}, we selected $Z_1 = \mu_1 X_1^{-1}$ as an initial value, where $X:\euclid^{\binom{q+1}{2}(n+m)}\rightarrow \symmetric^{q(n+m)}$ is a function
generated by combining constraints of problem~\cref{problem:psf_shfted} into a single positive semidefinite constraint, that is,
\begin{align}
    &X(A_1, \ldots, A_m, B_1, \ldots, B_n)\nonumber \\
    &\coloneqq \allowbreak (A_1+ r I_q)\oplus\cdots\oplus (A_m+ r I_q)\oplus (B_1+ r I_q)\cdots\oplus (B_n+ r I_q),
    \label{constraint:oplus}
\end{align}
where $\oplus$ denotes the direct sum of matrices.
\par
Turning to the assumptions we made, the above parameter rules satisfy~\cref{assumption:compleconv}. The validity of other assumptions is less clear.
\subsection{Results}
\subsubsection{Starting from the neighborhood of the saddle point}
\label{subsubsection:saddle_neighbor}
We set $(m,n,q,r) = (5, 5, 4, 0.3)$ in problem~\eqref{problem:psf_shfted} and see the behavior of both algorithms in their first 300 iterations.
\Cref{tab:num_exp} summarizes the results of the experiments. ``N/A'' in the second column means not applicable because the algorithm of the corresponding row is not equipped with negative curvature.
All of the results indicate that the algorithm with negative curvature obtains \tteal{the lower function value}.
Even though the algorithm with negative curvature takes more time to reach 300 iterations, less function value is also confirmed if we compare the two algorithms with respect to function value vs. time.
\begin{table}[htbp]
    \centering
    \caption{Decreasing-NC-PDIPM with negative curvature (the bottom row of each cell) vs. without negative curvature (the top row of each cell).}
    \label{tab:num_exp}
    \begin{tabular}{l|lll}
        Seed & \# negative curvature & Function value & Time              \\ \hline
        1    & N/A                   & \num{1.40e1}  & \SI{1.90}{\second} \\
             & 3                  & \num{4.34e-1}  & \SI{2.10}{\second} \\ \hline
        2    & N/A                   & \num{5.43e-1}  & \SI{1.80}{\second} \\
             & 7                    & \num{2.00e-1}  & \SI{2.17}{\second} \\ \hline
        3    & N/A                   & \num{5.95e-1}  & \SI{3.80}{\second} \\
             & 5                     & \num{6.17e-5}  & \SI{4.19}{\second} \\ \hline
        4    & N/A                   & \num{6.75e1}     & \SI{1.94}{\second} \\
             & 5                    & \num{1.64e-1}  & \SI{2.29}{\second} \\ \hline
        5    & N/A                   & \num{1.69e1}  & \SI{1.80}{\second} \\
            & 5                     & \num{3.26e-1}  & \SI{2.43}{\second} \\ \hline
        6    & N/A                   & \num{1.29e1}  & \SI{1.92}{\second} \\
            &  5                   & \num{6.32e-2}  & \SI{2.24}{\second}
    \end{tabular}
\end{table}
\subsubsection{Comparison with projected gradient methods}
\label{subsubsection:comparison_with_pgd}
 We also implement projected gradient methods (PGD) for comparison in numerical experiments. Each iteration of PGD is defined by
\begin{equation*}
    \xkp = \Proj_S \paren{\xk - \alpha \nabla f(\xk)},
\end{equation*}
    where $\alpha >0$, $S$, and $\Proj_S$  denote the stepsize, the feasible region, and the orthogonal projection operator onto $S$, that is,
\begin{equation*}
    \Proj_S (t) = \argmin_{s \in S} \norm{s-t}^2,
\end{equation*}
respectively. We set the stepsize $\alpha = 10^{-3}$.
\par 
We compute the projection for problem~\eqref{problem:psf_shfted} in the following way.
As for projection onto $\symmetric^q_{+}$, it can be calculated efficiently. Let $T\in \symmetric^q$, which can be decomposed as $T= UDU^\top$, where $D$ is a diagonal matrix and $U$ is an orthogonal matrix. We define $D_{+}$ by
\begin{equation*}
    D_{+} = \diag \paren{\max \sbra{0, D_{11}}, \ldots ,\max\sbra{0, D_{qq}}}.
\end{equation*}
 Then, projection of $T\in \symmetric^q$ onto $\symmetric^q_{+}$ is written as follows~\cite[pp. 399]{BoydVandenberghe2004}:
\begin{equation*}
    \Proj_{\symmetric^q_{+}}(T) = \argmin_{S\in \symmetric^q_{+}}
    \norm{S-T} = U D_{+} U^\top.
\end{equation*}
In particular, when $Y \in \symmetric^{q(m+n)}$ is written as
\begin{equation*}
    Y = 
    T_1 \oplus \cdots \oplus T_{m+m}
\end{equation*}
where $T_{i} \in \symmetric^q$ for all $i \in \sbra{1,\ldots, m+n}$, the projection of $Y$ onto $\symmetric^{q(m+n)}_{+}$ is 
\begin{equation*}
    \Proj_{\symmetric^{q(m+n)}_{+}}(Y) = 
      \Proj_{ \symmetric^q_{+}} (T_1) \oplus \cdots \oplus \Proj_{ \symmetric^q_{+}}(T_{m+n}).
\end{equation*}
Thus, when we apply PGD to problem~\eqref{problem:psf_shfted}, we can calculate projection of $X(A_1, \ldots, A_m, \allowbreak B_1, \ldots, B_n)$, which is defined in~\eqref{constraint:oplus}, to $\symmetric^{q(m+n)}_{+}$ by a series of eigenvalue decomposition of $q\times q$ symmetric matrices.
\begin{figure}[htbp]
    \centering
    \includegraphics[width=\linewidth]{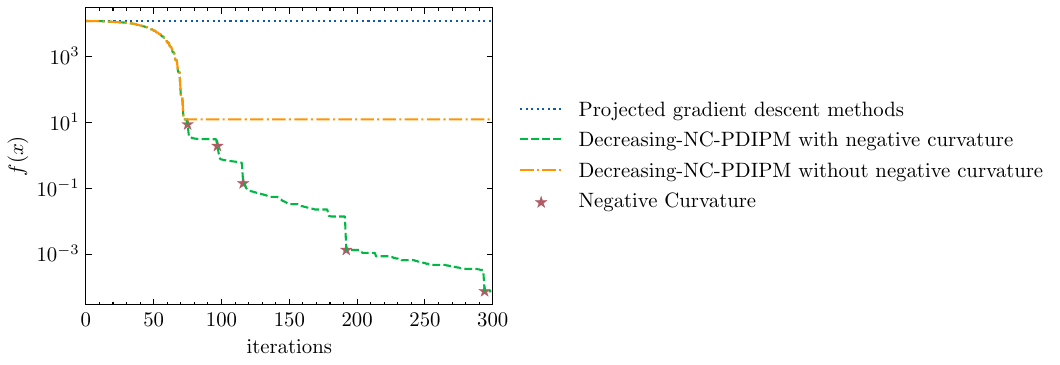}
    \caption{Comparison with PGD}
    \label{figure:comparison}
\end{figure}
\par
\Cref{figure:comparison}  shows that Decreasing-NC-PDIPM with negative curvature could successfully escape from the strict saddle point although PGD could not escape from it and Decreasing-NC-PDIPM without negative curvature could escape from it but 
become stuck at a high function value $f(x)$. The usage of negative curvature described as \textcolor{negative_curvature}{$\boldsymbol{\star}$} decreases the function value significantly.

\section{Discussion and conclusion}
\label{section:conclu}
We have proposed a primal-dual IPM with convergence to SOSPs of NSDPs using directions of negative curvature. Utilizing ``local'' Lipschitz continuities of derivatives of $\pmunu$ makes it possible to establish a worst-case iteration complexity of our primal-dual IPM for computing {\ASOSP} of NSDPs. The difficulty of our analysis comes from the sigma-term that appears in the second-order necessary conditions~\eqref{eq:sigma-term} for NSDPs. Our numerical experiments show that the use of negative curvature is beneficial in nonconvex settings.
\par
We believe that our primal-dual IPM can be extended to the equality-constrained NSDPs and it is possible to use trust-region methods in place of negative curvature methods in inner iterations.

\paragraph{\rm \bf Data Availability Statement }
The source code utilized in the numerical experiments can be accessed at \url{https://github.com/Mathematical-Informatics-5th-Lab/Decreasing-NC-PDIPM}. There is no conflict of interest in writing the paper.
\begin{acknowledgements}
We are very grateful to the two anonymous reviewers and the associate editor for their valuable comments and suggestions. In addition, we also thank Dr. Hiroshi Yamashita and Professor Hiroshi Yabe for helpful discussions.
\end{acknowledgements}

\bibliographystyle{spmpsci}      
\bibliography{abbreviated.bib}   
\appendix
\setcounter{section}{0}
\renewcommand{\thesection}{\Alph{section}}
\setcounter{equation}{0}
\renewcommand{\theequation}{\thesection-\arabic{equation}}

\section{Proofs of \texorpdfstring{\cref{lemmabreak:localLipschitz}}{}}
\label{appendix:local_lipschitz}
\subsection{Proof of~\texorpdfstring{\cref{lemma:invX}}{}}
\begin{proof}
    Notice that
    \[
        \fnorm{\invXell - \invXx}  \allowbreak=\allowbreak\fnorm{\invXell\paren{\Xell - \Xx}\invXx}            \allowbreak
        \leq\allowbreak \fnorm{\invXell}\allowbreak\fnorm{\invXx}
        \allowbreak\fnorm{\Xell-\Xx}.
    \]
    Since $\fnorm{\invXx}\le \fnorm{\invXell}$ follows from \cref{lemma:feasiblity}, inequality\,\eqref{eq:Xlip} yields
    $\bigl\| \invXell -\invXx\bigr\|_{\mathrm{F}}
        \allowbreak \leq 2 L_0 \fnorm{\invXell}^2
        \norm{\xell-x}.
    $ The proof is complete.
    \qed
\end{proof}
\subsection{Proof of~\texorpdfstring{\cref{lemma:nabalxLips}}{}}
\begin{proof}
    Using~\cref{eq:gradfx}, we obtain
    \begin{align*}
         & \norm{\nabla_x \psimunu{\xell, \Zell} - \nabla_x \psimunu{x, \Zell}}\leq \underbrace{\norm{\nabla f(\xell)-\nabla f(x)}}_{(A)} \\
         & \;+\underbrace{\nu\norm{\paren{\Aast (\xell)\Zell - \Aast (x)\Zell}}}_{(B)}
        +\underbrace{\mupn\norm{\Aast (\xell)\invXell - \Aast (x)\invXx}}_{(C)}
        .
    \end{align*}
    In what follows, we evaluate upper bounds of $(A),(B),$ and $(C)$.\newline
    {\it Evaluation of $(A):$}
    From \eqref{eq:gradf}, we have
    \begin{equation}
        (A) \leq L_1 \norm{\xell - x}.
        \label{ap2:last1}
    \end{equation}
    \newline
    {\it Evaluation of $(B):$}
    \begin{equation}
        \begin{aligned}
            (B)
             & \leq
            \nu \sum_{i=1}^n\abs{\trace{\Zell \paren{\Acali\paren{\xell}-\Acali\paren{x}}}}                                               \\
             & \leq \nu \sum_{i=1}^n \fnorm{\Zell}\fnorm{\Acali\paren{\xell}-\Acali\paren{x}}\leq \nu L_1 \fnorm{\Zell} \norm{\xell - x},
        \end{aligned}
        \label{eq:l2B}
    \end{equation}
    where the third inequality follows from \eqref{eq:gradxlips}.
    \newline
    {\it Evaluation of $(C):$}
    Let
\begin{align*}
    (D) &\coloneqq \mupn \norm{\paren{\Aast (\xell)-\Aast (x)}\invXell},\\ 
  (E) &\coloneqq \mupn \allowbreak \norm{\Aast (x)\paren{\invXell-\invXx}}.
 \end{align*}
    As $(C)\leq (D) + (E)$, we will evaluate upper bounds of $(D)$ and $(E)$. With regard to $(D)$, we obtain from~\cref{eq:gradxlips} that
    \begin{equation*}
        \begin{aligned}
            \label{eq:l2D}
            (D)
             & \leq  \mupn \fnorm{\invXell}\sum_{i=1}^n \fnorm{\Aast (\xell)- \Aast(x)} \\
             & \leq \mupn L_1 \fnorm{\invXell}\norm{\xell - x}.
        \end{aligned}
    \end{equation*}
    Similarly, it follows from~\cref{eq:pdvXxiLips}, that
    \begin{align}
        \label{eq:l2E}
        (E) & \leq \mupn  \sum_{i=1}^n
        \abs{\trace{\Acali\paren{\xell}\paren{\invXell -\invXx}}}
        \nonumber                                                                                                                       \\
            & \leq  \mupn  \sum_{i=1}^n \fnorm{\Acali\paren{x}} \fnorm{\invXell -\invXx}
        \nonumber                                                                                                                       \\
            & \leq 2 \mupn L_0\sum_{i=1}^n \fnorm{\Acali \paren{x}} \fnorm{\invXell}^2 \norm{\xell - x}\quad \tag{by \cref{lemma:invX}} \\
            & \leq 2 \mupn L_0^2 \fnorm{\invXell}^2 \norm{\xell - x}.\quad \tag{by \cref{eq:hessbounded}}
    \end{align}
    Therefore
    \begin{equation}
        (C) \leq \paren{\mupn L_1 \fnorm{\invXell} + 2 \mupn L_0^2 \fnorm{\invXell}^2}\norm{\xell - x}.
        \label{eq:l2Cend}
    \end{equation}
    Lastly, from~\cref{ap2:last1,eq:l2B,eq:l2Cend}, we obtain
    \begin{equation}
        \begin{aligned}
             & \norm{\nabla_x \psimunu{\xell, \Zell}
                -\nabla_x \psimunu{x, \Zell}}\leq  \\
             & \quad 
            \paren{L_1 + \nu L_1\fnorm{\Zell}+ 2\mupn L_0^2\fnorm{\invXell}^2
                + \mupn L_1 \fnorm{\invXell}} \norm{\xell - x}.\nonumber
        \end{aligned}
    \end{equation}
    The proof is complete.
    \qed
\end{proof}
\subsection{Proof of~\texorpdfstring{\cref{lemma:nabalZLips}}{}}
\begin{proof}
    Note that
    \begin{align*}
         \fnorm{\nabla_Z \psimunu{\xell, \Zell} - \nabla_Z \psimunu{\xell, Z}} &= \mu \nu \fnorm{Z^{-1} - \Zell^{-1}} \\
         & \leq \mu \nu \fnorm{\Zell^{-1}}\fnorm{\Zell - Z} \fnorm{Z^{-1}},
    \end{align*}
    where the equality follows from \eqref{eq:gradfZ}.
    Since  $\fnorm{Z^{-1}} \leq 2\fnorm{\Zell^{-1}}$ from \cref{lemma:Zpd}, we have
    \[\fnorm{\nabla_Z \psimunu{\xell, \Zell}\allowbreak - \nabla_Z \psimunu{\xell, Z}} \allowbreak \leq 2\mu\nu \fnorm{\Zell^{-1}}^2 \fnorm{\Zell- Z}.\]
    The proof is complete.
    \qed
\end{proof}
\subsection{Proof of~\texorpdfstring{\cref{lemma:hessLips}}{}}
\begin{proof}
From \cref{eq:hessofpsi}, we have $$
        (\nabla^2_{xx} \psimunu{\xell, \Zell})_{ij} -
        (\nabla^2_{xx}\psimunu{x, \Zell})_{ij} \allowbreak {=} (A_{ij}) + (B_{ij}) + (C_{ij}),$$
    where
    \begin{align*}
        (A_{ij}) & :=\paren{\nabla^2 f(\xell)}_{ij}-\paren{\nabla^2 f\bigl(x\bigr)}_{ij},        \\
        (B_{ij}) & :=\mupn\Bigl( \trace{\Acali(\xell) \Xell^{-1}\Acalj(\xell) \Xell^{-1}}
        \\
                 & \quad -\trace{\Acali(x)\invXx\Acalj(x)\invXx}\Bigr),                            \\
        (C_{ij}) & := \trace{\paren{\paren{1+\nu} \mu \invXx- \nu \Zell}\pdv{X}{x_i}{x_j}\paren{x}} \\
                 & \quad
        -\trace{\paren{\paren{1+\nu} \mu \Xell^{-1}- \nu \Zell}\pdv{X}{x_i}{x_j}\paren{\xell}}.
    \end{align*}
    In what follows, we evaluate $(A_{ij})$, $(B_{ij})$, and $(C_{ij})$.\newline
    {\it Evaluation of $(A_{ij})$}: From~\cref{eq:hessf}, we have
    \begin{equation}
        \norm{\nabla^2 f(\xell) - \nabla^2 f(x)} \leq L_2 \norm{\xell-x}.
        \label{eq:evaluationofAij}
    \end{equation}
    \newline
    {\it Evaluation of $(B_{ij})$}: $(B_{ij})$ can be written as
        {\small
            \begin{align*}
                (B_{ij}) & =\underbrace{\mupn\paren{ \trace{\Acali(\xell) \invXell\Acalj(\xell)\invXell}-
                        \trace{\Acali(x) \invXx\Acalj(x)\invXell}}}_{(D_{ij})}
                \\
                         &\quad  +
                \underbrace{
                    \mupn\paren{
                        \trace{\Acali(x) \invXx \Acalj(x) \invXell}
                        -\trace{\Acali(x)\invXx\Acalj(x)\invXx}}}_{(E_{ij})}.
            \end{align*}}
    Here, $(D_{ij})$ can be bounded as
    \begin{align*}
         \abs{(D_{ij})} 
         & \leq \mupn
        \fnorm{\Xell^{-1}}
        \fnorm{\Acali(\xell) \Xell^{-1}\Acalj(\xell)-\Acali(x) \Xx^{-1}\Acalj(x)}
        \nonumber         \\
         & \leq
        \underbrace{\mupn
            \fnorm{\Xell^{-1}}
            \fnorm{\Acali(\xell) \Xell^{-1}\Acalj(\xell)-\Acali(x) \Xell^{-1}\Acalj(x)}}_{(F_{ij})}
        \\
         & \quad +
        \underbrace{
            \mupn
            \fnorm{\Xell^{-1}}
            \fnorm{\Acali(x) \Xell^{-1}\Acalj(x)-\Acali(x) \Xx^{-1}\Acalj(x)}}_{(G_{ij})}.
    \end{align*}
    In the above, $(F_{ij})$ can be bounded as
    \begin{align}
        \begin{split}
            (F_{ij})
            &\leq \mupn
            \fnorm{\Xell^{-1}}
            \fnorm{\Acali(\xell) \Xell^{-1}\Acalj(\xell)-\Acali(x) \Xell^{-1}\Acalj(\xell)}\\
            &\quad +\mupn
            \fnorm{\Xell^{-1}}
            \fnorm{\Acali(x) \Xell^{-1}\Acalj(\xell)-\Acali(x) \Xell^{-1}\Acalj(x)}\\
            &\leq  \mupn
            \fnorm{\Xell^{-1}}
            \fnorm{\Acali(\xell)  - \Acali(x) }
            \fnorm{ \Xell^{-1}\Acalj(\xell)}\\
            &\quad +
            \mupn
            \fnorm{\Xell^{-1}}\fnorm{\Acali(x)\Xell^{-1}}
            \fnorm{\Acalj(\xell)-\Acalj(x)}.
            \nonumber
        \end{split}
    \end{align}
    Taking the sum of $(F_{ij})$ over $i,j$ together with \cref{eq:pdvXxiLips,eq:gradxlips} yields
    \begin{align}
        \begin{split}
            \sum_{i=1}^n \sum_{j=1}^n (F_{ij})
            &\leq 2\mupn L_0 L_1{\fnorm{\invXell}^2}
            \norm{\xell -x}.
            \label{eq:17211end}
        \end{split}
    \end{align}
    For $(G_{ij})$, we obtain, from~\cref{lemma:invX},
    \begin{align}
        \begin{split}
            (G_{ij})&\leq \mupn
            \fnorm{\Xell^{-1}}\fnorm{\Acali(x)}\fnorm{\Acalj(x)}
            \fnorm{ \Xell^{-1}- \Xx^{-1}}
            \\
            &\leq
            2\mupn L_0
            \fnorm{\Xell^{-1}}^3\fnorm{\Acali(x)}\fnorm{\Acalj(x)}
            \norm{\xell-x},
            \nonumber
        \end{split}
    \end{align}
    which together with~\cref{eq:pdvXxiLips} yields
    \begin{align}
        \begin{split}
            \sum_{i=1}^n\sum_{j=1}^n
            (G_{ij})\leq
            2\mupn L_0^3
            \fnorm{\invXell}^3
            \norm{\xell-x}.
            \label{eq:l7212end}
        \end{split}
    \end{align}
    By combining inequalities~\eqref{eq:17211end} and~\eqref{eq:l7212end}, the sum of $(D_{ij})$ is bounded from above by
    \begin{align}
        \begin{split}
            &\sum_{i=1}^n \sum_{j=1}^n
            \abs{(D_{ij})}\\
            &\leq
            2\paren{
                \mupn L_0^3
                \fnorm{\invXell}^3+
                \mupn L_0 L_1{\fnorm{\invXell}^2}
            }\norm{\xell-x}.
            \label{eq:1721end}
        \end{split}
    \end{align}
    With regard to $(E_{ij})$, we have
    \begin{equation*}
        \abs{\paren{E_{ij}}}
        \leq \mupn\fnorm{\Acali(x)}\fnorm{\Acalj(x)}\fnorm{\Xx^{-1}}
        \fnorm{ \Xell^{-1}-  \Xx^{-1}}.
    \end{equation*} By taking the sum over $i,j$ of $\abs{\paren{E_{ij}}}$ together with~\cref{eq:pdvXxiLips}, ~\cref{lemma:invX}, and~\cref{lemma:feasiblity}, we have
    \begin{align}
        \begin{split}
            \sum_{i=1}^n \sum_{j=1}^n
            \abs{\paren{E_{ij}}}
            &\leq 4 \mupn
            L_0^3 \fnorm{\Xell^{-1}}^3
            \norm{\xell - x}.
            \label{eq:1722end}
        \end{split}
    \end{align}
    Combining inequalities~\cref{eq:1721end} and~\cref{eq:1722end} yields
    \begin{align}
        \begin{split}
            &\sum_{i=1}^n \sum_{j=1}^n \abs{\paren{B_{ij}}}\\
            &\leq \paren{6\mupn L_0^3
                \fnorm{\invXell}^3+ 2\mupn
                L_0 L_1 {\fnorm{\Xell^{-1}}^2}}
            \norm{\xell - x}.
            \label{eq:l72end}
        \end{split}
    \end{align}
    \newline {\it Evaluation of $(C_{ij})$}: Note that $(C_{ij})$ can be written as
    \begin{align*}
        \paren{C_{ij}}
         & =\underbrace{\nu\trace{\Zell \paren{\pdv{X}{x_i}{x_j}\paren{\xell}-\pdv{X}{x_i}{x_j}\paren{x}}}}_{(J_{ij})}     \\
         & \quad +\underbrace{\mupn \trace{\invXx\pdv{X}{x_i}{x_j}\paren{x}-\invXell\pdv{X}{x_i}{x_j}\paren{\xell}}}_{(K_{ij})}.
    \end{align*}
    By the Cauchy-Schwartz inequality, we have
    \begin{align}
        \begin{split}
            \sum_{i=1}^n \sum_{j=1}^n\abs{\paren{J_{ij}}}
            &\leq \nu \sum_{i=1}^n \sum_{j=1}^n\fnorm{\Zell}\fnorm{\pdv{X}{x_i}{x_j}\paren{\xell}-\pdv{X}{x_i}{x_j}\paren{x}}
            \\
            &\leq \nu L_2 \fnorm{\Zell} \norm{\xell -x}.
            \label{eq:l731end}
        \end{split}
    \end{align}
    With regard to the sum of $\abs{(K_{ij})}$,
    \begin{align}
         & \sum_{i=1}^n \sum_{j=1}^n
        \abs{(K_{ij})}\nonumber                                                                                 \\
         & =
        \sum_{i=1}^n \sum_{j=1}^n
        \abs{\mupn  \trace{\paren{\invXx-\invXell}\pdv{X}{x_i}{x_j}\paren{x}}}\nonumber
        \\
         &\quad +
        \sum_{i=1}^n \sum_{j=1}^n\abs{
        \mupn \trace{\invXell\paren{\pdv{X}{x_i}{x_j}\paren{x}-\pdv{X}{x_i}{x_j}\paren{\xell}}}}
        \nonumber \\
         & \leq 2\mupn  L_0 L_1 {\fnorm{\invXell}^2}
        \norm{\xell -x} +\mupn L_2 \fnorm{\invXell}\norm{\xell -x},
        \label{eq:l732end}
    \end{align}
    where the inequality follows from~\cref{eq:hessbounded},~\cref{as:hessXlips}, and~\cref{lemma:invX}. Moreover, combining~\cref{eq:l731end,eq:l732end} yields
    \begin{equation}
        \begin{aligned}
             & \sum_{i=1}^n \sum_{j=1}^n
            \abs{\paren{C_{ij}}}                                                                                   \\
             & \leq  \paren{   \nu L_2 \fnorm{\Zell}  + 2\mupn  L_0 L_1{\fnorm{\invXell}^2}+\mupn L_2 \fnorm{\invXell
                }}\norm{\xell -x}
            \label{eq:l73end}.
        \end{aligned}
    \end{equation}
    Lastly, {
       we have 
       \begin{align*}
             &\norm{\nabla^2_{xx} \psimunu{\xell, \Zell} - \nabla^2_{xx} \psimunu{x, \Zell}}\\
         &\le \norm{\nabla^2 f(\xell)-\nabla^2 f(x)} + \|(B_{ij})_{i,j}\|+ \|(C_{ij})_{i,j}\|\\
         &\le \norm{\nabla^2 f(\xell)-\nabla^2 f(x)} + \|(B_{ij})_{i,j}\|_{\rm F}+ \|(C_{ij})_{i,j}\|_{\rm F}\\
         &\le \norm{\nabla^2 f(\xell)-\nabla^2 f(x)} + \sum_{i=1}^n\sum_{j=1}^n\left|(B_{ij})\right|+\sum_{i=1}^n\sum_{j=1}^n\left|(C_{ij})\right|,
     \end{align*}
which together with \cref{eq:evaluationofAij,eq:l72end,eq:l73end} implies}
    \begin{align*}
         & \norm{\nabla^2_{xx} \psimunu{\xell, \Zell} - \nabla^2_{xx} \psimunu{x, \Zell}} \nonumber \\
         & \leq
        \biggl( L_2+\nu L_2 \fnorm{\Zell}  \\
         & \quad \quad  +  
        \mupn\paren{ L_2 \fnorm{\invXell} +
            4 L_1 L_0\fnorm{\invXell}^2 + 6 L_0^3 \fnorm{\invXell}^3}\biggr)\\
        & \quad \cdot  \norm{\xell - x}
    \end{align*}
    The proof is complete.
    \qed
\end{proof}

\end{document}